\documentclass[11pt]{article}

\title{The Morse index theorem for regular Lagrangian systems}
\author{Chaofeng Zhu\\ \\Department of Mathematics,
Massachusetts Institute of Technology
\\Cambridge, MA 02139-4307 USA \\and \\
Nankai Institute of Mathematics,
Nankai University\\Tianjin 300071, The People's Republic of China}

\date{}

\begin{document}

\maketitle

\newtheorem{Def}{Definition}[section]
\newtheorem{Th}{Theorem}[section]
\newtheorem{Prop}{Proposition}[section]
\newtheorem{Not}{Notation}[section]
\newtheorem{Lemma}{Lemma}[section]
\newtheorem{Rem}{Remark}[section]
\newtheorem{Cor}{Corollary}[section]

\def\s{\section}
\def\ss{\subsection}

\def\d{\begin{Def}}
\def\t{\begin{Th}}
\def\p{\begin{Prop}}
\def\n{\begin{Not}}
\def\la{\begin{Lemma}}
\def\r{\begin{Rem}}
\def\c{\begin{Cor}}
\def\ee{\begin{equation}}
\def\aa{\begin{eqnarray}}
\def\y{\begin{eqnarray*}}
\def\bd{\begin{description}}

\newcommand{\proof}[1]{\noindent {\em Proof.}$\quad$ {#1} $\hfill\Box$
                                \vspace{2ex}}
\def\ed{\end{Def}}
\def\et{\end{Th}}
\def\ep{\end{Prop}}
\def\en{\end{Not}}
\def\el{\end{Lemma}}
\def\er{\end{Rem}}
\def\ec{\end{Cor}}
\def\eee{\end{equation}}
\def\eaa{\end{eqnarray}}
\def\ey{\end{eqnarray*}}
\def\ebd{\end{description}}

\def\nn{\nonumber}
\def\bp{{\bf Proof.}\hspace{2mm}}
\def\qe{\hfill{\rm Q.E.D.}}
\def\lj{\langle}
\def\rj{\rangle}
\def\dd{\diamond}
\def\ox{\mbox{}}
\def\lb{\label}

\def\K{{\bf K}}
\def\R{{\bf R}}
\def\C{{\bf C}}
\def\Z{{\bf Z}}
\def\N{{\bf N}}
\def\Ua{{\bf U}}
\def\Da{{\bf D}}

\def\B{{\cal B}}
\def\Ca{{\cal C}}
\def\A{{\cal A}}
\def\Ga{{\cal G}}
\def\Ha{{\cal H}}
\def\Ka{{\cal K}}
\def\Pa{{\cal P}}
\def\J{{\cal J}}
\def\g{{\bf g}}

\def\gl{{\rm gl}}
\def\Gl{{\rm Gl}}
\def\Sp{{\rm Sp}}
\def\sp{{\rm sp}}
\def\U{{\rm U}}
\def\O{{\rm O}}
\def\G{{\rm G}}
\def\H{{\rm H}}
\def\P{{\rm P}}
\def\D{{\rm D}}
\def\T{{\rm T}}
\def\Sa{{\rm S}}
\def\rel{\;{\rm rel.}\;}
\def\vp{\epsilon}
\def\mod{\;{\rm mod}\;}
\def\Lie{{\rm Lie}}
\def\diag{{\rm diag}}
\def\im{{\rm im}\;}
\def\Lag{{\rm Lag}}
\def\Gr{{\rm Gr}}
\def\span{{\rm span}}
\def\Stab{{\rm Stab}\;}
\def\sign{{\rm sign}}
\def\sf{{\rm sf}}
\def\ind{{\rm ind}}
\def\rank{{\rm rank}}
\def\Sg{{\Sp(2n,\C)}}
\def\Na{{\cal N}}
\def\dist{{\rm dist}}

\begin{abstract}
{\it In this paper, we prove a Morse index theorem for the index form of
regular Lagrangian system with selfadjoint boundary condition. }

\end{abstract}

\s{Introduction}

Let $(M,g)$ be an $n$-dimensional Riemannian manifold. The classical Morse
Index Theorem states that the number of conjugate points along a geodesic
$\gamma:[a,b]\to M$ counted with multiplicities is equal to the index of
the second variation of the Riemannian action functional
$E(c)=\frac{1}{2}\int_a^bg(\dot c,\dot c)dt$
at the critical point $\gamma$, where $\dot c$ denotes $\frac{d}{dt}c$.
Such second variation is called the index form
for $E$ at $\gamma$. The theorem has later been extended in several
directions (see \cite{AgSa,Am,Du,PiT1,PiT2,Sm,Uh} for versions of this theorem
in different contexts). In \cite{Du} of 1976, J. J. Duistermaat proved his general Morse
index theorm for positive definite Lagrangian system with selfadjoint boundary conditions.
In \cite{AgSa} of 1996, A. A. Agrachev and A. V. Sarychev studied the Morse index
and rigidity of the abnormal sub-Riemannian geodesics. In \cite{PiT1, PiT2} of 2000,
P. Piccione and D. V. Tausk proved a version of the Morse index theorem
for geodesics in semi-Riemannian geodesics with both endpoints varies on two submanifolds
of $M$ under some nondegenerate conditions (cf. Theorems 5.2 and 5.9 in \cite{PiT2}).
In this paper, we will prove a general version of Morse index theorem
for regular Lagrangian system with selfadjoint boundary conditions, and show how the
indices varies under different choices of the frames. The relation between the indices
for two different boundary conditions is a easy corollary of Proposition \ref{p01m2.2}
below. In order to get such a general
theorem for regular Lagrangian system, one has to overcome the
following difficulties.
\bd
\item[(1)] The multiplicities of conjugate points may be meaningless.
\item[(2)] The Morse index of the index form may be infinity.
\item[(3)] The corresponding second order operators may have different domains.
\ebd
We overcome these difficulties by using the notions of Maslov-type indices
(see \cite{A,CLM,Lo7, LZh}) and the spectral flow
(see \cite{APS1,BoFu98,DZ2,RS2}).

Let $M$ be a smooth manifold of dimension $n$, points in its tangent bundle
$TM$ will be denoted by $(m,v)$, with $m\in M$, $v\in T_mM$. Let $f$ be a
real-valued $C^3$ function on an open subset $Z$ of $\R\times TM$. Then

\ee \lb{eq01m1.6}
E(c)=\int_0^Tf\left(t,c(t),{\dot c}(t)\right)dt
\eee
defines a real-valued $C^2$ function $E$ on the space of curves

\ee \lb{eq0m1.7}
{\cal C}=\left\{c\in C^1([0,T],M);(t,c(t),{\dot c}(t))\in Z
\;\mbox{\rm for}\;\mbox{\rm all}\; t\in[0,T]\right\}.
\eee
The set ${\cal C}$ is a $C^2$ Banach manifold modeled on the Banach space
$C^1([0,T],\R^n)$ with its usual topology of uniform convergence of
the curves and their derivatives.

Boundary conditions will be introduced by restricting $E$ to the set of curves

\ee \lb{eq01m1.8}
{\cal C}_N=\{c\in{\cal C};(c(0),c(T))\in N\},
\eee
where $N$ is a given smooth submanifold of $M\times M$. The most
familiar example are $N=\{(m(0),m(T)\}$ and
$N=\{(m_1,m_2)\in M\times M;m_1=m_2\}$.
In the general case ${\cal C}_N$ is a smooth submanifold of ${\cal C}$ with
tangent space equal to
\ee \lb{eq01m1.9}
T_c{\cal C}_N=\left\{\delta c\in C^1([0,T],c^*TM);
(\delta c(0),\delta c(T))\in T_{(c(0),c(T))}N\right\}.
\eee
$c\in {\cal C}_N$ is called {\bf stationary curve} for the boundary
condition $N$ if the restriction of $E$ to ${\cal C}_N$ has a stationary
point at $c$, that is,
if $DE(c)(\delta c)=0$ for all $\delta c\in T_c{\cal C}_N$. For such a curve
$c$ is of class $C^2$.

Let $c$ is of class $C^2$. Then the second order differential $D^2E(c)$ of
$E$ at $c$ is symmetric bilinear form on $T_c{\cal C}_N$, which is called the
{\bf index form} of $E$ at $c$ with respect to the boundary condition $N$.
In general the Morse index of $D^2E(c)$ will be infinite. In order to get
a well-defined integer, we introduce the following concept.

Assume that
$f$ is a {\bf regular Lagrangian}, that is,
\ee \lb{eq01m1.10}
D_v^2f(t,m,v)\;\mbox{{\rm is}}\;\mbox{{\rm nondegenerate}}
\;\mbox{{\rm for}}\;\mbox{{\rm all}}\;(t,m,v)\in Z.
\eee
Here $D_v$ denotes differential of functions on $Z$ with respect to
$v\in T_mM$, keeping $t$ and $m$ fixed.
The condition (\ref{eq01m1.10}) is called the {\bf Legendre condition}.

Let $H=H^1(T_c{\cal C}_N)$ be the $H^1$ completion of $T_c{\cal C}_N$. By
Sobolev embedding theorem, $H\subset C([0,T],c^*TM)$. Then
$D^2E(c)$ is well defined on $H$. In local coordinates, we have

\aa
D^2E(c)(X,Y)&=&\int_0^T((D_v^2f({\tilde c}(t)){\dot\alpha}
+D_mD_v({\tilde c}(t))\alpha,{\dot\beta})\nn\\
& &+(D_vD_m({\tilde c}(t)){\dot\alpha},\beta)
+(D_mD_m({\tilde c}(t))\alpha,\beta))dt,
\lb{eq01m1.11}\eaa
where $X,Y\in H$, $\alpha$, $\beta$ are the local coordinate expression of $X$, $Y$
defined by
$X=(\alpha,\partial m)$, $Y=(\beta,\partial m)$, $\partial m$ is
the natural frame of $T_mM$,
and we have use the abbreviation

\ee \lb{eq01m1.12}
{\tilde c}(t)=(t,c(t),\dot c(t)).
\eee

In general $\partial m$ and $\alpha$ is not globally well-defined along $c$.
Choosing a $C^1$ frame $e$ of $T_c{\cal C}_N$. Such a frame can be
obtained by the parallel transformation of the induced connection
on $c^*TM$ of a connection on $TM$ (for example, the Levi-Civita
connection with respect to a semi-Riemannian metric on $TM$). Then
in local coordinates, there is a $C^1$ path $a(t)\in \gl(n,\R)$
which is nondegenerate for all $t$ such that $\partial
m=(a(t),e(t))=a(t)^*e(t)$, where $a(t)^*$ denotes the transpose
conjugate of $a(t)$. The vector fields $X,Y\in H$ along $c$ can be
written as $X=(x,e)$, $Y=(y,e)$, where $x,y\in H^1([0,T],\R^n)$
and $(x(0),x(T)), (y(0), y(T))\in R$, $R$ is defined by
\footnote{In this paper, all vectors are viewed as column vectors. For a pair of vectors
$x,y\in\C^n$, $(x,y)$ has two meanings: one is the standard
Hermitian inner product of $x,y$, the other is the the vector
$(x^*,y^*)^*\in\C^{2n}$. The readers can easily see it from the
content.}
$$R=\{(x,y)\in\R^{2n};((x,e(0)),(y,e(T)))\in T_{(c(0),c(T))}N\}.$$
So we have

\ee \lb{eq01m1.13}
x=a\alpha,\quad \dot x=a{\dot\alpha}+{\dot a}\alpha,\quad
y=a\beta,\quad \dot y=a{\dot\beta}+{\dot a}\beta.
\eee
Substitute (\ref{eq01m1.13}) to (\ref{eq01m1.11}), we get the
following form of the the index form:

\ee \lb{eq01m1.14}
D^2E(c)(X,Y)=\int_0^T\left((p\dot x+qx,\dot y)+(q^*\dot x,y)
+(r x,y)\right)d t,
\eee
where $p,q,r\in C([0,T],\gl(n,\R))$, $p$ is of class $C^1$,
$p(t)=p(t)^*$, $r(t)=r(t)^*$, $p(t)$ are invertible
for all $t\in[0,T]$, and $^*$ denotes
the conjugate transpose. Now define
\ee \lb{eq01m1.15}
{\cal I}_s(x,y)=\int_0^T\left((p\dot x+sqx,\dot y)+(sq^*\dot x,y)
+(s r x,y)\right)d t,\quad s\in[0,1],
\eee
where $x,y\in H^1([0,T],R^n)$ and $(x(0),x(T)), (y(0),y(T))\in R$.
Since $p$ is of class $C^1$ and is nondegenerate,
the relative Morse index $I({\cal I}_0,{\cal I}_1)\equiv-\sf\{{\cal I}_s\}$
is a well-defined finite integer. When $p$ is positive definite,
$I({\cal I}_0,{\cal I}_1)$ is the Morse index of $D^2E(c)$
(where we only require that $p$ is continuous).

The main results in this paper is the following.

Let $p,q,r\in C([0,1]\times[0,T],\gl(n,\C))$ be families of matrices
such that $p$ is of class $C^1$, $p_s(t)=p_s(t)^*$, $r_s(t)=r_s(t)^*$,
and $p_s(t)$ is invertible for all $s\in[0,1]$ and $t\in[0,T]$.

Let $R$ be a subspaces of $\C^{2n}$. Let $H_R$ be the Hilbert space
defined by

\ee \lb{eq01m1.3}
H_R=\left\{x\in H^1([0,T],\C^n);(x(0), x(T))\in R\right\}.
\eee
For each $s\in[0,1]$, let ${\cal I}_s$ be the index form defined by

\ee \lb{eq01m1.4}
{\cal I}_s(x,y)=\int_0^T\left((p_s\dot x+q_sx,\dot y)+(q_s^*\dot x,y)
+(r_sx,y)\right)dt,\qquad x,y\in H_R.
\eee

Let $J=\pmatrix{0&-I_n\cr I_n&0\cr}$, $I_n$ is the
identity matrix on $\R^n$. When there is no confusion we will omit the subindex
of identity matrices. Let

\ee \lb{eq01m1.1}
b_s(t)=\pmatrix{p_s^{-1}(t)&-p_s^{-1}(t)q_s(t)\cr
-q_s^*(t)p_s^{-1}(t)&q_s^*(t)p_s^{-1}(t)q_s(t)-r_s(t)\cr}.
\eee
For each $s\in[0,1]$, let $\gamma_s(t)$ be the fundamental solution of the
linear Hamiltonian equation

\ee \lb{eq01m1.2}
\dot u=Jb_s(t)u.
\eee

Define

\y R^b&=&\{(x,y)\in\C^{2n};(x,-y)\in R^{\perp}\},\\
W(R)&=&\{(x,y,z,u)\in\C^{4n};(x,-z)\in R^{\perp},(y,u)\in R\}.
\ey

Then we have

\t \lb{th01m1.1} Let $\sf\{{\cal I}_s;0\le s\le 1\}$ be the
spectral flow of ${\cal I}_s$,
and $i_{W(R)}(\gamma)$ be the Maslov-type index of $\gamma$
defined below. Then we have

\ee \lb{eq01m1.5}
-\sf\{{\cal I}_s;0\le s\le 1\}=i_{W(R)}(\gamma_1)-i_{W(R)}(\gamma_0).
\eee
\et

Assume that $q_0(t)=r_0(t)=0$ for all $t\in[0,T]$. Then we have
$b_0(t)=\diag(p_0^{-1}(t),0)$ and $\gamma_0(t)=\pmatrix{I&0\cr
\int_0^tp_0^{-1}(s)ds&I\cr}$ for all $t\in[0,T]$.

\t \lb{th01m1.2} Let $P\in C([0,T],\gl(n,\C))$ be a path of selfadjoint
matrices. Define
$\gamma(t)=\pmatrix{I&0\cr P(t)&I\cr}$ for all $t\in[0,T]$.
Then we have

\ee \lb{eq01m1.18}
i_{W(R)}(\gamma)=m^+(P(T)|_S)-m^+(P(0)|_S),
\eee
where $m^+$ denotes the Morse positive index, and
$$S=\{x\in\C^n;(x,x)\in\R^b\}.$$\et

As a special case, we get the following theorem of J. J.
Duistermaat \cite{Du}.

\c \lb{c01m1.1} Assume that $p_1(t)$ are positive definite for all
$t\in[0,T]$. Then we have

\ee \lb{eq01m1.17}
m^-({\cal I}_1)=i_{W(R)}(\gamma_1)-\dim_{\C}S,
\eee

where $m^-$ denotes the Morse (negative) index, and
$$S=\{x\in\C^n;(x,x)\in\R^b\}.$$ \ec

Let $a(t)\in \Gl(n,\C)$, and
$$R^{'}=\{(x,y)\in\C^{2n}; (a(0)x,a(T)y)\in R\}.$$
After change of frame $x\longmapsto ax$, $I_1(ax,ay)$ defines a quadratic
form on $H_{R^{'}}$ and we can get the corresponding $p_1^{'},q_1^{'}$ and
$r_1^{'}$. By (\ref{eq01m1.1}) and (\ref{eq01m1.2}) we get the corresponding
$b_1^{'}$ and $\gamma_1^{'}$.

\t \lb{th01m1.3} We have the following
\ee \lb{eq01m1.16}
i_{W(R^{'})}(\gamma_1^{'})- i_{W(R)}(\gamma_1)
=\dim_{\C}(\Gr(I)\cap R^{'})-\dim_{\C}(\Gr(I)\cap R),
\eee
where $\Gr(I)$ denotes the graph of $I$.\et

The paper is organized as follows. In {\S}2, we discuss the
properties of the spectral flow. In {\S}3, we discuss the
properties of the Maslov-type indices. In {\S}4, we prove our main
results.

{\bf Acknowledgements}
This work was done when the author visited Professor Tian Gang at MIT. The
author is most grateful to Professor Liu Chun-gen, Professor Tian Gang and Professor
Zhang Weiping for simulating discussions and comments, and MIT for the
hospitality and nice research air there.

\s{Spectral flow}

\ss{Definition of the spectral flow}

The spectral flow for a one parameter family of linear selfadjoint Fredholm
operators is introduced by Atiyah-Patodi-Singer \cite{APS1} in their study of
index theory on manifolds with boundary. Since then other significant
applications have been found. In \cite{DZ2} the notion of the spectral flow was
generalized to the higher dimensional case by X. Dai and W. Zhang. In
\cite{Zh,ZhL} it is generalized to more general operators.

Let $X$ be a Banach space. We denote the set of closed operators, bounded
linear operators and compact linear operators on $X$ by $\Ca(X)$, $\B(X)$
and ${\cal CL}(X)$ respectively. For $A\in\Ca(X)$, an operator $B$
is called {\bf $A$-compact} if $\Da(A)\subset\Da(B)$, and view $B$ as operator
from $\Da(A)$ to $X$, is compact, where $\Da(A)$ is the domain of $A$ with
the graph norm of $A$.

According to Atiyah-Patodi-Singer \cite{APS1}, we define

\d (cf. Definition 1.3.6 of \cite{Zh} and Definition 2.6 of \cite{ZhL})
\lb{d01m2.1} (1) Let $l$ be a real dimension $1$ cooriented embedded $C^1$
submanifold of $\C$ without boundary. Let $A$ be in $\Ca(X)$. $A$ is said to
be {\bf admissible} with respect to $l$ if the spectrum of $A$ near $l$
lies on a compact subset of $l$ and is of finite algebraic multiplicity.
If $\infty$ is a limit point of $l$, we require $A\in\B(X)$.
Let $P_l(A)$ be the spectral projection of $A$ on $l$. The {\bf nullity}
$\nu_{l}(A)$ of $A$ with respect to $l$ is defined to be
$\nu_{l}(A)=\dim_{\C}\im P_l(A)$.
All such $A$ will be denoted by $\A_l(X)$.

(2) Let $A_s$, $0\le s\le 1$ be a curve in $\A_l(X)$. The {\bf spectral slow}
$\sf_l\{A_s \}$ of $A_s$ counts the algebraic multiplicities of
the spectral of $A_s$ cross the manifold $l$, i.e., the number of
the spectral lines of $A_s$ cross $l$ from the negative part of $\C$ near $l$ to
the non-negative part of $\C$ near $l$ minus the number of the
spectral lines of $A_s$ cross $l$ from the non-negative part of
$\C$ near $l$ to the negative part of $\C$ near $l$.
When $l=\sqrt{-1}(-K,K)$ and $A_s\in\A_l(X)$ be such that
$\sigma(A_s)\cap \sqrt{-1}\R\subset\sqrt{-1}(-K,K)$ for all $s$,
where $(-K,K)$ ($K>0$) denotes the open interval on $\R$ and
$\sigma(A)$ denotes the spectrum of $A$,
we set $\sf\{A_s\}=\sf_l\{A_s\}$.

(3) Let $l=\sqrt{-1}\R$, $A\in\A_l(X)$ and $B\in {\cal CL}(X)$,
or $X$ is a Hilbert space, $A$ is a selfadjoint Fredholm operator
with compact resolvent and $B$ is bounded selfadjoint operator,
The {\bf relative Morse index} of the pair $A$, $A+B$ is defined by

\ee \lb{eq01m2.10}
I(A,A+B)=-\sf\{A+sB\}.
\eee

(4) When $l=\sqrt{-1}\R$ and $A\in\A_l(X)$,
or $A$ is selfadjoint Fredholm on Hilbert space
$X$ and $l=\sqrt{-1}(-\vp,\vp)$,
we call the algebraic multiplicity of the spectrum of $A$ in the right side and
the left side $l$ the {\bf Morse positive index} and the {\bf Morse
(negative) index}, and denote them by $m^+(A)$ and $m^-(A)$.
The {\bf signature} $\sign(A)$ is defined by
$\sign(A)=m^+(A)-m^-(A)$.
\ed

\ss{Calculation of the spectral flow}

Let $X$ be a complex Banach space, $\gamma$ be a $C^1$ curve in $\C$ which
bounds a bounded open subset $\Omega$ of $\C$. Let $A_s$, $s\in(-\vp,\vp)$,
where $\vp>0$, be a curve in ${\cal C}(X)$. Assume that
$\gamma\cap\sigma(A_s)=\emptyset$ for all $s\in(-\vp,\vp)$, where
$\sigma(A_s)$ denotes the spectral of $A_s$. Set $A_0=A$, and

\ee \lb{eq01m2.21}
P_s\equiv P_{\gamma}(A_s)
=-\frac{1}{2\pi \sqrt{-1}}\int_{\gamma}R(\zeta,A_s)d\zeta,
\eee
where $R(\zeta,A_s)=(A_s-\zeta I)^{-1}$, $\zeta\in\C\setminus\sigma(A_s)$
is the resolvent of $A_s$. Then $P_s^2=P_s$. Set $P_0=P$. Assume that
$\im P\subset \Da(A_s)$, for all $s\in(-\vp,\vp)$,
$\im P$ is a finitely dimensional subspace of $X$,
and $\frac{d}{ds}|_{s=0}(A_sP)=B$. Then $B$ is
bounded. Let $f$ be a polynomial. Then $P_sf(A_s)P_s$ is uniformly bounded on
any compact subsets of $(-\vp,\vp)$, and

\ee \lb{eq01m2.22}
P_sf(A_s)P_s=-\frac{1}{2\pi \sqrt{-1}}
\int_{\gamma}f(\zeta)R(\zeta,A_s)d\zeta.\eee
Set $R_s=(I-(P_s-P)^2)^{-\frac{1}{2}}$. Then $R_sP=PR_s$ and $R_sP_s=P_sR_s$.
Set

\y
U_s^{'}&=&P_sP+(I-P_s)(I-P), \qquad U_s=U_s^{'}R_s,\\
V_s^{'}&=&PP_s+(I-P)(I-P_s), \qquad V_s=V_s^{'}R_s.
\ey
Then we have

\y
U_sV_s&=&V_sU_s=I, \\
U_sP&=&P_sU_s=P_sR_sP,\\
PV_s&=&V_sP_s=PR_sP_s.\ey

\la \lb{l01m2.11} We have

\ee \lb{eq01m2.23}
\frac{d}{ds}|_{s=0}(U_s^{-1}P_sA_sP_sU_s)
=\frac{1}{2\pi \sqrt{-1}}
\int_{\gamma}R(\zeta,A)BR(\zeta,A)d\zeta.\eee
If $(PAP)(PB)=(PB)(PAP)$, then we have

\ee \lb{eq01m2.24}
\frac{d}{ds}|_{s=0}(U_s^{-1}P_sA_sP_sU_s)=PB.\eee
\el

\bp By the definition of $U_s$ and $V_s$ we have

$$U_s^{-1}P_sA_sP_sU_s=V_sP_sP_sU_s=PR_sP_sA_sP_sR_sP.$$
By (\ref{eq01m2.22}) we have

\ee \lb{eq01m2.25}
(P_sf(A_s)P_s-Pf(A)P)P=\frac{1}{2\pi \sqrt{-1}}
\int_{\gamma}f(\zeta)R(\zeta,A)(A_sP-AP)R(\zeta,A)d\zeta.\eee
Since $A_s$, $s\in(-\vp,\vp)$ is a curve in ${\cal C}(X)$, and $\im P$
has finite dimension, we have

\ee \lb{eq01m2.26}
\frac{d}{ds}|_{s=0}(P_sf(A_s)P_sP)=\frac{1}{2\pi \sqrt{-1}}
\int_{\gamma}f(\zeta)R(\zeta,A)BR(\zeta,A)d\zeta.\eee

Take $f=1$, we have $\frac{d}{ds}|_{s=0}P_sP$ exists. By the definition of
$R_s$ we have $\frac{d}{ds}|_{s=0}R_sP=0$. Now formal calculation shows

$$\frac{d}{ds}|_{s=0}(PR_sP_sA_sP_sR_sP)
=\frac{1}{2\pi \sqrt{-1}}
\int_{\gamma}PR(\zeta,A)BR(\zeta,A)d\zeta.$$
The assumption that $\im P$ has finite dimension shows the calculation is
right.

When $(PAP)(PB)=(PB)(PAP)$, we have

\y
\frac{d}{ds}|_{s=0}(PR_sP_sA_sP_sR_sP)
&=&\frac{1}{2\pi \sqrt{-1}}
\int_{\gamma}PR(\zeta,A)R(\zeta,A)Bd\zeta\\
&=&P^2B\\
&=& PB.\ey
\qe

By the above Lemma \ref{l01m2.11}, the proof of Theorem 4,1 in \cite{ZhL}
also works for the following more general proposition.

\p \lb{p01m2.1}
Let $l$ be a $C^1$ submanifold of $\R$ without boundary.
Let $A_s$, $-\vp\le s\le \vp$ ($\vp>0$),be a curve in $\A_l(X)$.
Set $P=P_l(A_0)$ and $A=A_0$ Assume that $\im P\subset \Da(A_s)$ and
$B=\frac{d}{ds}\mid_{s=0}A_sP$ exists.
Assume that
\begin{equation}\lb{eq01m2.1}
(PAP)(PB)-(PB)(PAP)=0,
\end{equation}
where $PAP, PB\in\B(\im P)$, and $PB:\im P\to\im P$ is hyperbolic.
Then there is a $\delta\in(0,\vp)$ such that
$\nu_l(A_s)=0$ for all $s\in [-\delta,0)\cup (0,\delta]$ and
\begin{eqnarray}
\lb{eq01m2.2}
& &\sf_l\{A_s,0\le s\le\delta\}=-m^{-}(PB:\im P\to \im P),\\
\lb{eq01m2.3}
& &\sf_l\{A_s,-\delta\le s\le 0\}=m^+(PB:\im P\to \im P).
\end{eqnarray}
\ep\qe

Now we consider two special cases.

\la \lb{l01m2.6} Let $X$ be a Hilbert space. Let $A_s, -\vp\le s\le \vp$ be a
curve of selfadjoint Fredholm operators with constant domain $D$
such that $A_s\le A_t$ for all $-\vp<s<t<\vp$. Assume that
$\frac{d}{ds}A_sx$ exist for all $s\in(-\vp,\vp)$, $x\in D$, or $A_s$ is
bounded for all $s\in(-\vp,\vp)$.
Then for $s<0$ with $|s|$ small, $\dim_{\C}\ker A_s$
is constant and $A_s$ has no positive small eigenvalue. For $s>0$ small,
$\dim_{\C}\ker A_s$ is constant and $A_s$ has no negative
eigenvalue whose absolute value is small.
\el

\bp Firstly assume that $\frac{d}{ds}A_sx$ exist for all
$s\in(-\vp,\vp)$, $x\in D$.
Let $\lambda_1(s)\le\ldots\le\lambda_k(s)$ be the spectral lines of $A_s$
for $|s|$ small such that $\lambda_1(0)=\ldots=\lambda_k(0)=0$.
Fix $j=1,\ldots k$ and $t$ with $|t|$ small.
Pick $x(s)\in\ker(A_s-\lambda_j(s)I)$ such that
$\|x_s\|=1$. Then every subsequence of $x_s$, $s\to t$ has a convergent
subsequence. Let $s_n$ be the subsequence of $s$, $s\to t$ such that
$$a_j(t)\equiv\liminf_{s\to t}\frac{\lambda_j(s)-\lambda_j(t)}{s-t}
=\lim_{n\to\infty}\frac{\lambda_j(s_n)-\lambda_j(t)}{s_n-t}$$
and $\lim_{n\to\infty}x(s_n)\to x$. Then $x\in\ker(A_t-\lambda_j(t)I)$.
So we have

\y a_j(t)&=&\lim_{n\to\infty}\frac{\lambda_j(s_n)-\lambda_j(t)}{s_n-t}\\
&=&\lim_{n\to\infty}\frac{\lambda_j(s_n)-\lambda_j(t)}{s_n-t}(x(s_n),x)\\
&=&\lim_{n\to\infty}\left(\frac{(A_{s_n}-A_t)x}{s_n-t}x,x(s_n)\right)\\
&=&\left(\frac{d}{ds}|_{s=t}A_sx,x\right)\\
&\ge& 0.
\ey
Hence $\lambda_j(s)\le \lambda_j(t)$ for $s<t$ and $|s|,|t|$ small, and our
results follows.

Now assume that $A_s$ is bounded. For $t>0$ small, consider the curve
$A_0+(t-s)(A_t-A_s)$, $0\le s\le t$. By the above arguments, $A_t$ has
no negative eigenvalue near $0$.
\qe

By Lemmas \ref{l01m2.1} below, \ref{l01m2.6} and the definition of the
spectral flow we have

\p \lb{p01m2.3} Let $X$ be a Hilbert space.

(1) Let $A_s, 0\le s\le 1$ be a
curve of selfadjoint Fredholm operators with constant domain $D$
such that $A_s\le A_t$ for all $0\le s<t\le 1$. Assume that
$\frac{d}{ds}A_sx$ exist for all $s\in[0,1]$, $x\in D$, or $A_s$ is
bounded for all $s\in[0,1]$. Then we have
\ee \lb{eq01m2.15}
\sf\{A_s\}=\sum_{0<s\le 1}\left(\dim_{\C}\ker A_s
-\lim_{t\to s^-}\dim_{\C}\ker A_t\right)\ge 0.
\eee

(2) Let $A\in{\cal C}(X)$ be a selfadjoint operator with compact resolvant and
$B_s\in\B(X)$, $0\le s\le 1$
be a curve of selfadjoint operators such that $B_s\le B_t$ for all
$0\le s<t\le 1$. Then we have

\ee \lb{eq01m2.16}
\sf\{A+B_s\}=\sum_{0<s\le 1}\left(\dim_{\C}\ker (A+B_s)
-\lim_{t\to s^-}\dim_{\C}\ker (A+B_t)\right)\ge 0.
\eee
\ep
\qe

Similarly we have

\p \lb{p01m2.11} Let $X$ be a Hilbert space and $l=(1-\vp,1+\vp)$
($\vp>0$ small). Let $A_s\in\B(X)$, $0\le s\le 1$ be a curve of
unitary operators. Assume that $A_s-I$ is Fredholm and
$\sqrt{-1}A_s^{-1}{ \dot A}_s\le 0$ for all
$s\in[0,1]$. Then we have
\ee \lb{eqm012.27}
\sf_l\{A_s\}=\sum_{0<s\le 1}\left(\dim_{\C}\ker (A_s-I)
-\lim_{t\to s^-}\dim_{\C}\ker (A_t-I)\right)\ge 0.
\eee
\ep
\qe

\ss{Spectral flow for curves of quadratic forms}

The spectral flow for curves of selfadjoint Fredholm operators has some
special properties as finite dimensional case.

The following lemma is Corollary 2.2 of \cite{ZhL}.

\la \lb{l01m2.1}
Let $X$ be a Hilbert space. Let $A$ be a selfadjoint Fredholm
operator on $X$ with compact resolvent, and $B$ be a bounded selfadjoint
operator on $X$. Set $K=(\|A\|+I)^{-1}$. Then we have
\begin{equation}\lb{eq01m2.4}
I(AK,AK+KB\}=I(A,A+B),
\end{equation}
where $AK$, $AK+KB$ are linear operators defined on the Hilbert space
$V=D(|A|^{\frac{1}{2}})$ with graph norm
$$\|x\|_V=(\| |A|^{\frac{1}{2}}x\|^2_{X}+\|x\|^2_{X})^{\frac{1}{2}}.$$
\el\qe

\la \lb{l01m2.2} Let $X$ be a Hilbert space. Let $A_s$, $0\le s\le 1$
be a curve of closed selfadjoint Fredholm operators.
Then for any curve $P_s\in\B(X)$ of invertible operators, we have
\ee\lb{eq01m2.5}
\sf\{P_sP_s^*A_s\}=\sf\{P_s^*A_sP_s\}=\sf\{A_s\}.
\eee
\el

\bp Since $A_s$ is a curve of closed selfadjoint Fredholm operators and $P_s$
is a curve of bounded invertible operators, the family
$P_s^*A_sP_s$, $0\le s\le 1$, is a curve of closed selfadjoint
Fredholm operators.
By the definition of the spectral flow we have
\aa
\sf\{P_sP_s^*A_s\}&=&\sf\{P_s(P_s^*A_sP_s)P_s^{-1}\}\nn\\
&=&\sf\{P_s^*A_sP_s\}.
\lb{eq01m2.6}\eaa

Since $P_s^*A_tP_s$ are selfadjoint Fredholm operators and
$\dim_{\C}\ker (P_s^*A_tP_s)=\dim_{\C}\ker A_t$, we have
\aa
\sf\{P_s^*A_sP_s\}&=&\sf\{P_0^*A_sP_0\}+\sf\{P_s^*A_1P_s\}\nn\\
&=&\sf\{P_0^*A_sP_0\}\nn\\
&=&\sf\{P_1^*A_sP_1\}.
\lb{eq01m2.7}\eaa

Let $Q_s$ be curves of bounded positive operators on $X$ satisfying $Q_0=I$,
$Q_1^2=P_0P_0^*$. By (\ref{eq01m2.6}) and (\ref{eq01m2.7}) we have

\y \sf\{P_s^*A_sP_s\}&=& \sf\{P_0^*A_sP_0\}\\
&=&\sf\{P_0P_0^*A_s\}\\
&=&\sf\{Q_1A_sQ_1\}\\
&=&\sf\{Q_0A_sQ_0\}\\
&=&\sf\{A_s\}.\ey
\qe

The above lemma leads the following definition.

\d \lb{d01m2.3}
Let $X$ be a Hilbert space. Let ${\cal I}_s$, $0\le s\le 1$ be a family of
quadratic forms. Assume that ${\cal I}_s(x,y)=(A_sx,y)$ for all $x,y\in X$,
where $A_s$ is a curve of bounded selfadjoint Fredholm operators.

(1) The {\bf spectral flow} $\sf\{{\cal I}_s\}$ of ${\cal I}_s$ is
defined to be the spectral flow $\sf\{A_s\}$.

(2) If $A_1-A_0\in{\cal CL}(X)$, the {\bf relative Morse index}
$I({\cal I}_0,{\cal I}_1)$ is defined to be the
relative Morse index $I(A_0,A_1)$.\ed

Based on this observation we can prove the following lemma.

\la \lb{l01m2.3}
Let $X$ be a Hilbert space. Let $A_s\in\B(X)$, $0\le s\le 1$ be a curve
of selfadjoint Fredholm operators and ${\cal I}_s$ be quadratic forms
defined by ${\cal I}_s(x,y)=(A_sx,y)$ for all $x,y\in X$. Assume that
$P_s\in\B(X)$, $0\le s\le 1$ is a curve of operators such that
$P_s^2=P_s$ and ${\cal I}_s(x,y)=0$
for all $x\in\im P_s$, $y\in\im Q_s$, where $Q_s=I-P_s$.

\ee \lb{eq01m2.8}
\sf\{{\cal I}_s\}=\sf\{{\cal I}_s|_{\im P_s}\}+\sf\{{\cal I}_s|_{\im Q_s}\}.
\eee
\el

\bp Let $R_s=P_s^*P_s+Q_s^*Q_s$. Since
$P_s+Q_s=I$ and $P_s^2=I$, we have

$$R_s=\frac{I}{2}+2(\frac{I}{2}-P_s^*)(\frac{I}{2}-P_s)>0.$$
So $R_s^{-1}A_s$ are Fredholm operators. Now consider the new inner product
$(R_sX,y)$, $x,y\in X$ on $X$. For this inner product $P_s$ is an orthogonal
projection. Let $B_s\in \B(\im P)$ and $C_s\in\B(\im Q_s)$ satisfy
${\cal I}_s(x,y)=(B_sx,y)$ for all $x,y\in \im P_s$
and ${\cal I}_s(x,y)=(C_sx,y)$ for all $x,y\in \im Q_s$ respectively.
By the fact that
$\im P_s$ and $\im Q_s$ are ${\cal I}_s$-orthogonal, we have

$$(A_sx,y)=(R_s(B_s\oplus C_s)x,y),\quad\forall x,y\in X.$$
So $R_s^{-1}A_s=B_s\oplus C_s$, and $B_s$, $C_s$ are Fredholm operators.
By Lemma \ref{l01m2.2}
and the definition of the spectral flow we have

\y
\sf\{{\cal I}_s\}&=&\sf\{R_s^{-1}A_s\}\\
&=&\sf\{B_s\}+\sf\{C_s\}\\
&=&\sf\{{\cal I}_s|_{\im P_s}\}+\sf\{{\cal I}_s|_{\im Q_s}\}.
\ey
\qe

\r \lb{r01m2.1} Here we allow the Hilbert space $\im P_s$
continuous varying. By Lemma I.4.10 in \cite{Ka}, for $t\in[0,1]$
being close enough to $s$, there are
invertible operators $U_{s,t}\in\B(X)$  such that

$$P_tU_{s,t}=U_{s,t}P_s,\qquad
U_{s,t}\to I, \quad {\rm as}\; t\to s.
$$
So locally we can define the spectral flow of $B_t$ as that of
$U_{s,t}^{-1}B_tU_{s,t}:\im P_s\to\im P_s$ ($s$ fixed), and globally
add them up.
\er

\la \lb{l01m2.4}
Let $X$ be a Hilbert space, and $M$ be a closed subspace
with finite codimension. Let $A\in\B(M)$ be a selfadjoint Fredholm operator
and ${\cal I}(x,y)=(Ax,y)$ for all $x,y\in M$. Let $N_0$ and $N_1$ be
subspace of $X$ such that $X=M\oplus N_0=M\oplus N_1$.
Define ${\cal I}_k$ on $X$, $k=0,1$  by

$${\cal I}_k(x+u,y+v)=(Ax,y),\quad\forall x,y\in M,\forall u,v\in N_k.$$
Then we have $I({\cal I}_0,{\cal I}_1)=0$.
\el

\bp Without loss of generality, we assume that $N_0$ is the orthogonal
complement of $M$. Set $A_0=\diag(A,0)$.  Let $B:N_1\to N_0$ be a linear
isomorphism. Define $P_1\in \B(X)$ by $P_1(x+y)=x+By$ for all $x\in M$,
$y\in N_1$. Set $A_1=P_1^*A_0P_1$. Then $P_1$ is invertible, $P_1-I$ is compact,
and ${\cal I}_k(x,y)=(A_kx,y)$ for all $x,y\in X$ and $k=0,1$.
Let $P_s\in\B(X)$,
$0\le s\le 1$ be a curve of invertible operators such that $P_0=I$ and
$P_s-I$ are compact. By the definition of the spectral flow we have

\y I({\cal I}_0,{\cal I}_1)&=&I(A_0,A_1)\\
&=&-\sf\{P_s^*A_0P_s)\\
&=&0.
\ey
\qe

The following proposition gives a generalization of Proposition 5.3 in
\cite{AgSa} and a formula of M. Morse.

\p \lb{p01m2.2}
Let $X$ be a Hilbert space and $A\in\B(X)$ be a selfadjoint Fredholm
operator. Let $P$ be an orthogonal projection such that $\ker P$ is
of finite dimensional. Let ${\cal I}$ be quadratic form on $X$ defined by
${\cal I}(x,y)=(Ax,y)$, $x,y\in X$. Set $M=\im P$ and $N$ be the
${\cal I}$-orthogonal complement of $M$: $N=\{x\in X;{\cal I}(x,y)=0,
\forall y\in M\}$. Then we have

\ee \lb{eq01m2.9}
I(PAP,A)=m^-({\cal I}|_N)
+\dim_{\C}\ker {\cal I}|_N-\dim_{\C}\ker {\cal I}.
\eee
\ep

\bp Since $PAP-A$ is a finite rank operator, $sPAP+(1-s)A$ are selfadjoint
Fredholm operators. We divide our proof into three steps.

{\bf Step 1.} Equation (\ref{eq01m2.9}) holds when
$\ker A=0$ and $N\subset M$.

In this case, set $M_0=\ker {\cal I}|_M$, $M_1$ be the orthogonal
complement of $M_0$ in $M$, and $P_0$, $P_1$ be the orthogonal projection
onto $M_0$, $M_1$ respectively. Then $P_0$ is of finite rank and $P=P_0+P_1$.
Let $N_1$ be the ${\cal I}$-orthogonal complement of $M_1$. Since
${\cal I}|_{M_1}$ is nondegenrate, $M_1\cap N_1=\{0\}$. Moreover we have

\y
\dim_{\C}N_1&=&\dim_{\C}\ker (AP_1)-\ind (AP_1)\\
&=&\dim_{\C}\ker P_1-\ind A-\ind P_1\\
&=&\dim_{\C}\ker P_1<\infty.
\ey
So $X=M_1\oplus N_1$. By the fact that ${\cal I}$ is nondegenrate,
${\cal I}|_{N_1}$ is nondegenrate. Clearly $M_0\subset N_1$ and $M_0$ is the
${\cal I}|_{N_1}$-orthogonal complement of $M_0$. $N_1$ has an orthogonal
decomposition $N_1=N^+\oplus N^-$ such that $N^+$ and $N^-$ is
${\cal I}$-orthogonal, ${\cal I}|_{N^+}>0$ and ${\cal I}|_{N^-}<0$. Let
$P^{\pm}$ be the orthogonal projection onto $N^{\pm}$. Then $P^{\pm}|_{M_0}$
is isomorphism. So we have

$$\dim_{\C}N_1=2\dim_{\C}M_0=2m^-({\cal I}|_{N_1}).$$
Let ${\cal I}_1$ be defined by ${\cal I}_1(x+u,y+v)={\cal I}(x,y)$ for all
$x,y\in M_1$, $u,v\in N_1$.
By Lemma \ref{l01m2.3} and \ref{l01m2.4} we have

\y
I(PAP,A)&=&I(PAP, P_1AP_1)+I(P_1AP_1,A)\\
&=&I(P_1AP_1,A)\\
&=&I({\cal I}_1,{\cal I})\\
&=&I({\cal I}_1|_{M_1},{\cal I}|_{M_1})
+I({\cal I}_1|_{N_1},{\cal I}|_{N_1})\\
&=&m^-({\cal I}|_{N_1})\\
&=&\dim_{\C}\ker M_0.
\ey

{\bf Step 2.} Equation (\ref{eq01m2.9}) holds if $M+N=X$.

In this case we have

$$\ker {\cal I}|_N=\ker {\cal I}=M\cap N.$$

Firstly we assume that $\ker A=\{0\}$. Then $M\cap N=\{0\}$. Let
${\cal I}_1$ be defined by ${\cal I}_1(x+u,y+v)={\cal I}(x,y)$ for all
$x,y\in M$, $u,v\in N$. By Lemma \ref{l01m2.3} and \ref{l01m2.4} we have

\y
I(PAP,A)&=&I({\cal I}_1,{\cal I})\\
&=&I({\cal I}_1|_M),{\cal I}|_M)
+I({\cal I}_1|_N),{\cal I}|_N)\\
&=&m^-({\cal I}|_N).
\ey

In the general case, we apply the above special case by taking quotient space
with $\ker A$ and get $I(PAP,A)=m^-({\cal I}|_N)$.

{\bf Step3.} Equation (\ref{eq01m2.9}) holds.

Firstly we assume that $\ker A=\{0\}$. Let $Q$
be the orthogonal projection onto $M+N$. Then the ${\cal I}$-orthogonal
complement is $\ker_{\C}({\cal I}|_N)$. By the above two steps we have

\y
I(PAP,A)&=&I(PAP,QAQ)+I(QAQ,A)\\
&=&m^-({\cal I}|_N)
+\dim_{\C}\ker {\cal I}|_N.
\ey

In the general case, we apply the above special case by taking quotient space
with $\ker A$ and get (\ref{eq01m2.9}).
\qe

At the end of this subsection we gives the following formula which will
be used below.

\la \lb{l01m2.8}
Let $X$ be a Hilbert space and $A_s\in\Ca(X)$, $0\le s\le 1$ be a curve of
Fredholm operators. Let $H=X\oplus X$ and $B_s\in\Ca(H)$ be defined by
$B_s=\pmatrix{0&A_s^*\cr A_s&0\cr}$. Then we have

\ee \lb{eq01m2.11}
\sf\{B_s\}=\dim_{\C}\ker A_1-\dim_{\C}\ker A_0.
\eee
\el

\bp Note that for $\lambda\in\R$, $\lambda\in\sigma(B_s)$ if and only if
$\lambda^2\in\sigma(A^*A)$, and the algebraic multiplicity of them are the same if
$|\lambda|\ne 0$ is small. Moreover we have

\y
\dim_{\C}\ker B_s&=&\dim_{\C}\ker A_s+\dim_{\C}\ker A_s^*\\
\ind A_s=\ind A_0&=&\dim_{\C}\ker A_s-\dim_{\C}\ker A_s^*.
\ey
By the definition of the spectral flow we have

\y \sf\{B_s\}&=&\frac{1}{2}
\left(\dim_{\C}\ker B_1-\dim_{\C}\ker B_0)\right)\\
&=&\dim_{\C}\ker A_1-\dim_{\C}\ker A_0.
\ey
\qe

\s{Maslov-type index theory}

\ss{Introduction to Maslov index}

We begin with the definition of the Lagrangian Grassmannian of a
symplectic Hilbert space.

\d \lb{d01m3.1}
Let $X$ be a Hilbert space. Let $j\in\B(X)$ be an invertible skew
selfadjoint operator. Set $\omega(x,y)=(jx,y)$ for all $x,y\in X$.
The form $\omega$ is called the {\bf (strong) symplectic structure}
on $X$, and the space $(X, \omega)$ is called {\bf symplectic Hilbert space}.
The {\bf linear symplectic group} $\Sp(X,\omega)$ is defined to be
$$\Sp(X,\omega)=\{M\in \B(X);M^*jM=j\}.$$
Let $(X_l,\omega_l)$, $l=1,2$ be two symplectic Hilbert spaces. The space of
linear symplectic maps $\Sp(X_1,X_2)$ is defined to be
$$\Sp(X_1,X_2)=\{M\in\B(X_1,X_2);\omega_2(Mx,My)=\omega_1(x,y)\}.$$\ed

Set $A=(-j^2)^{\frac{1}{2}}$ and $J=A^{-1}j$. Then $(Ax,y)$, $x,y \in X$ is
an equivalent Hermitian metric on $X$ and $J$ is a complex structure on $X$
compatible with $\omega$, i.e., $J^2=-I$ and $\omega(x,Jy)=(Ax,y)$,
$x,y\in X$ is an equivalent Hermitian metric. All such $J$ forms a
contractable space. So we can replace the original metric on $X$ with
$A$.

\d \lb{d01m3.2}
Let $(X,\omega)$ be a symplectic Hilbert space.

(1) For any subspace $\Lambda$ of $X$, the {symplectic complement}
$\Lambda^{\omega}$ is defined to be
$$\Lambda^{\omega}=\{y\in X;\omega(x,y)=0,\forall x\in\Lambda\}.$$

(2) A subspace $\Lambda$ of $X$ is called {\bf Lagrange}
if $\Lambda^{\omega}=\Lambda$. The {\bf Lagrangian Grassmannian}
${\cal L}(X,\omega)$ consists of all the Lagrange subspaces of $(X,\omega)$.

(3) Let $\Lambda,\Lambda^{'}\in{\cal L}(X,\omega)$ be two Lagrange subspace.
The pair $(\Lambda,\Lambda^{'})$ is called a {\bf Fredholm pair} if
$\Lambda\cap\Lambda^{'}$ is of finite dimension and $\Lambda+\Lambda^{'}$ is
a finite
codimensional subspace of $X$. The {\bf Fredholm Lagrangian Grassmannian}
${\cal FL}_{\Lambda}(X,\omega)$ consists of all the Lagrange subspace
$\Lambda^{'}$ of $X$ such that $(\Lambda,\Lambda^{'})$ is a Fredholm pair.
\ed

The following lemma is well-known.

\la \lb{l01m3.1}
Let $(X,\omega)$ be a symplectic Hilbert space. Assume that there is
an invertible skew selfadjoint operator $J\in\B(X)$ such that $J^2=-I$ and
$\omega(x,y)=(Jx,y)$. Let $X_1=\ker (J-\sqrt{-1}I)$ and
$X_2=\ker (J+\sqrt{-1}I)$. Let $\Lambda_0,\Lambda$ are two subspaces of $X$.
Then we have
\bd
\item[(1)] $\Lambda\in{\cal L}(X,\omega)$ if and only if there is an linear
isometric $U\in\U(X_1,X_2)$ such that $\Lambda$ is the graph of $\Gr(U)$
of $U$, where $\U(X_1,X_2)$ denotes the set of linear isometric between
$X_1$ and $X_2$, $\U(X)$ denotes the unitary group of $X$ and
$\Gr(U)$ denotes the graph of $U$. Specially, ${\cal L}(X,\omega)\ne\emptyset$
if and only if $\U(X_1,X_2)\ne\emptyset$.

\item[(2)] Let $U, U^{'}\in\U(X_1,X_2)$. Set $\Lambda=\Gr(U)$,
$\Lambda^{'}=\Gr(U^{'})$. Then $(\Lambda, \Lambda^{'})$ is a Fredholm pair
if and only if $U-U^{'}$ is Fredholm.
\ebd\el
\qe

Following \cite{BoFu98} we give the following definition.

\d  \lb{d01m3.3}
Let $(X,\omega)$ be a symplectic Hilbert space. Assume that there is
an invertible skew selfadjoint operator $J\in\B(X)$ such that $J^2=-I$ and
$\omega(x,y)=(Jx,y)$. Let $X_1=\ker (J-\sqrt{-1}I)$ and
$X_2=\ker (J+\sqrt{-1}I)$. Let $(\Lambda(s),\Lambda^{'}(s))$, $a\le s\le b$
be a curve of Fredholm pairs of Lagrange subspaces of $X$ such that
$\Lambda(s)=\Gr(U(s))$, $\Lambda^{'}(s)=\Gr(U^{'}(s))$, where
$U(s),U^{'}(s)\in\U(X_1,X_2)$. Let $l=(1-\vp,1+\vp)\subset\C$, be a interval
on $\R$, where $\vp>0$ is small. The coorientation of $l$ is
defined to be from the down half complex plane to the up half
complex plane.
The {\bf Maslov index} $i(\Lambda,\Lambda^{'})$ is defined to be the spectral
flow $-\sf_l\{U(s)^{'-1}U(s)\}$. It is independent of the compatible complex
structure $J$.
\ed

To calculate the Maslov indices, the standard method is the crossing form
(cf. \cite{Du} and \cite{RS2}).

Let $\Lambda(s)$, $a\le s\le b$ be a curve of Lagrange subspaces of $X$.
Let $W$ be a fixed Lagrangian complement of $\Lambda(t)$.
For $v\in\Lambda(t)$ and $|s-t|$ small,
define $w(s)\in W$ by $v+w(s)\in\Lambda(s)$. The form
$$Q(\Lambda,t)\equiv Q(\Lambda,W,t)(u,v)=\frac{d}{ds}|_{s=t}\omega(u,w(s)),
\quad\forall u,v\in \Lambda(t) $$
is independent of the choice of $W$.
Let $(\Lambda(s),\Lambda^{'}(s))$, $a\le s\le b$
be a curve of Fredholm pairs of Lagrange subspaces of $X$.
For $t\in[a,b]$, the {\bf crossing form}
$\Gamma(\Lambda,\Lambda^{'},t)$ is defined on $\Lambda(t)\cap\Lambda^{'}(t)$ by
$$\Gamma(\Lambda,\Lambda^{'},t)(u,v)=Q(\Lambda,t)(u,v)-Q(\Lambda^{'},t)(u,v),
\quad\forall u,v\in\Lambda(t)\cap\Lambda^{'}(t).$$
A {\bf crossing} is a time $t\in[a,b]$ such that
$\Lambda(t)\cap\Lambda^{'}(t)\ne\{0\}$. A crossing is called {\bf regular}
if $\Gamma(\Lambda,\Lambda^{'},t)$ is nondegenerate. It is called {\bf simple}
if in addition $\Lambda(t)\cap\Lambda^{'}(t)$ is one dimensional.

Now let $(X,\omega)$ be a symplectic Hilbert space with $\omega(x,y)=(jx,y)$,
$x,y\in X$, $j\in\B(X)$ and $j^*=-j$. Then we have a symplectic Hilbert space
$(H=X\oplus X,(-\omega)\oplus \omega)$. For any $M\in\Sp(X,\omega)$, its graph
$\Gr(M)$ is a Lagrange subspace of $H$. The following lemma is
Lemma 3.1 in \cite{Du}.

\la \lb{l01m3.2}
Let $M(s)\in\Sp(X,\omega)$, $a\le s\le b$ be a curve of linear symplectic maps.
Assume that $M(s)$ is differentiable at $t\in[a,b]$. Set
$B_1(t)=-j{\dot M(t)}M(t)^{-1}$ and $B_2(t)=-jM(t)^{-1}{\dot M(t)}$.
Then $B_1(t)$, $B_2(t)$ are selfadjoint, $B_2(t)=M(t)^*B_1(t)M(t)$ and we have

\ee \lb{eq01m3.1}
Q(\Gr(M),t)((x,M(t)x),(y,M(t)y))=(B_2(t)x,y).
\eee
\el
\qe



\c \lb{c01m3.1}
Let $(X,\omega)$ be a symplectic Hilbert space and
$(\Lambda(t),\Lambda^{'}(t))$, $a\le t\le b$ be a $C^1$ curve of Fredholm pairs
of Lagrange subspaces of $X$ with only regular crossing. Then we have

\ee \lb{eq01m3.2}
i(\Lambda,\Lambda^{'})=m^+(\Gamma(\Lambda,\Lambda^{'},a))
-m^-(\Gamma(\Lambda,\Lambda^{'},b))+\sum_{a<t<b}\sign(\Gamma(\Lambda,\Lambda^{'},t)).
\eee
\ec

\bp Pick an invertible skew selfadjoint operator $J\in\B(X)$ such that $J^2=-I$
and $\omega(x,y)=(Jx,y)$. Let $X_1=\ker (J-\sqrt{-1}I)$ and
$X_2=\ker (J+\sqrt{-1}I)$. By Lemma \ref{l01m3.1} there are curves of
isometric $U(t)$, $U^{'}(t)$ in $\U(X_1,X_2)$ such that $\Lambda(t)=\Gr(U(t))$
and $\Lambda^{'}(t)=\Gr(U^{'}(t))$. Apply Lemma \ref{l01m3.2} for $X_1$
with $j=-\sqrt{-1}I$, for any $x,y\in\ker(U(t)-U^{'}(t))$ and $t\in[a,b]$
we have
\y
\frac{d}{ds}|_{s=t}(-\sqrt{-1}U^{'-1}Ux,y)
&=&(-jU^{'-1}{\dot U}^{'}U^{'-1}Ux,y)+(jU^{'-1}{\dot U}x,y)\\
&=&(-jU^{'-1}{\dot U}^{'}x,y)+(U^{'-1}UjU^{-1}{\dot U}x,U^{'-1}Uy)\\
&=&(-jU^{'-1}{\dot U}^{'}x,y)+(jU^{-1}{\dot U}x,y)\\
&=&-\Gamma(\Lambda,\Lambda^{'},t)((x,Ux),(y,Uy)).
\ey

By Proposition \ref{p01m2.1} we obtain (\ref{eq01m3.2}).
\qe

By Proposition \ref{p01m2.11}, Lemma \ref{l01m3.2} and the proof of Corollary
\ref{c01m3.1} we have

\c \lb{c01m3.13} Let $(X,\omega)$ be a symplectic Hilbert space and
$M(s)\in\Sp(X,\omega)$, $a\le s\le b$ be a $C^1$ curve of linear symplectic maps.
Assume that $-j{\dot M(t)}M(t)^{-1}$ is semi-positive definite.
Let $H=(X\oplus X, (-\omega)\oplus \omega)$ and $W$ be a Lagrange
subspace of $X$. Then we have
\ee \lb{eq01m3.30}
i_W(M(t))=\sum_{0<s\le 1}\left(\dim_{\C}(\Gr(M(s))\cap W)
-\lim_{t\to s^-}\dim_{\C}(\Gr(M(t))\cap W)\right)\ge 0.
\eee
\ec\qe

\c \lb{c01m3.10} ({\bf Symplectic invariance})
Let $(X_l,\omega_l)$, $l=1,2$ be two symplectic Hilbert spaces,
and $M(t)\in\Sp(X_1,X_2)$, $a\le t\le b$ be a curve of linear symplectic
maps, and
$(\Lambda(t),\Lambda^{'}(t))$, $a\le t\le b$ be a curve of Fredholm pairs
of Lagrange subspaces of $X_1$. Then we have
\ee \lb{eq01m3.20}
i(M\Lambda,M\Lambda^{'})=i(\Lambda,\Lambda^{'}).
\eee
\ec

\bp Firstly assume that the curves $M$, $\Lambda$, $\Lambda^{'}$
are differentiable and the pairs $(\Lambda,\Lambda^{'})$ have only
regular crossing.
By the definition of the crossing form, for any $t\in[a,b]$  we have
$$\Gamma(M\Lambda,M\Lambda^{'},t)(M(t)u,M(t)v)=\Gamma(\Lambda,\Lambda^{'},t)(u,v)
\quad\forall u,v\in\Lambda(t)\cap\Lambda^{'}(t).$$
By Corollary \ref{c01m3.1}, equation (\ref{eq01m3.20}) holds. For the
general case, we can make a small perturbation of the curves $M$, $\Lambda$,
$\Lambda^{'}$ with their endpoints fixed such that they satisfy the above
condition. Then our result follows from the homotopy invariance rel. endpoints of
the Maslov indices.
\qe

\ss{The spectral flow formula and Maslov-type indices}

Let $A$ be a closed densely defined symmetric operator on a Hilbert space $H$
with domain $D_m$. Let $D_M$ be the domain of $A^*$. Define the inner product
on $D_M$ by

$$\langle x,y\rangle=(x,y)+(A^*x,A^*y).$$
Then $D_M$ is a Hilbert space and $D_m$ is a closed subspace of $D_M$.
The orthogonal complement of $D_m$ in $D_M$ is $\ker(A^{*2}+I)$.
Define $X=D_M/D_m$ and let $\gamma: D_M\to X$ be the canonical map. The map
$\gamma$ is called the {\bf abstract trace map}. Define $\omega:X\times X\to\C$
by
$$\omega(\gamma(x),\gamma(y))=(A^*x,y)-(x,A^*y),\quad\forall x,y\in D_M.$$
Then $(X,\omega)$ is a symplectic Hilbert space.

The following proposition is Theorem 5.1 in \cite{BoFu98}.

\p \lb{p01m3.1} {\rm (Spectral flow formula)}
Let $A$ be a closed densely defined symmetric operator on a Hilbert space $H$
with domain $D_m$ and let $C_t,t\in[a,b]$ be bounded. We assume that

1. $A$ has a selfadjoint extension $A_D$ with compact resolvent, where $D$
is the domain of $A$.

2. there exists a positive constant $a$ such that
$D_m\cap\ker (A^*+C_t+s)=\{0\}$ for any $|s|<a$ and any $t\in [a,b]$.

Then we have

\ee \lb{eq01m3.3}
\sf\{A_D+C_t\}=-i(\gamma(D),\gamma(\ker(A^*+C_t))).
\eee
\ep

{\bf Sketch of the proof.} \hspace{2mm}
Here we only consider the case that $C_t$ is of class $C^1$. The condition shows
that $(\gamma(D),\{\gamma(\ker(A^*+C_t)\})$ is a $C^1$ Fredholm pairs of
Lagrange subspaces of $X$. Let
$x,y\in\ker(A^*+C_t)$. Pick a Lagrange complement $W$ of
$\gamma(\ker(A^*+C_t))$. Then for $s$ with $|s-t|$ small, $W$ is also a
Lagrange complement of $\gamma(\ker(A^*+C_s))$ and there is
$y_s\in\ker(A^*+C_s)$ such that $y_s\to y$ when $s\to t$.
So we have
\y
\omega(\gamma(x),\gamma(y_s-y))&=&(A^*x,y_s-y)-(x,A^*(y_s-y)\\
&=&(-C_tx,y_s-y)-(x,-C_sy_s+C_ty)\\
&=&((C_s-C_t)x,y_s).
\ey
Differential it with respect to $s$ at $s=t$, we get

\ee\lb{eq01m3.4}
Q(\gamma(\ker(A^*+C_s)),t)=Q_t{\dot C}_tQ_t,
\eee
where $Q_t$ is the orthogonal projection from $H$ onto $\ker(A^*+C_t)$.

By Theorem 4.22 in \cite{RS2}, we can choose $\delta\in(0,a)$
sufficiently small such that $\ker(A+C_t+sI)=\{0\}$ for $t=a$ or $b$ and
$s\in(0,\delta]$, and $A_D+C_t+\delta I$ has only regular crossing.
By Proposition \ref{p01m2.1}, Corollary \ref{c01m3.1} and (\ref{eq01m3.4})
we have

\y
\sf\{A_D+C_t\}&=&\sf\{A_D+C_t+\delta I\}\\
&=&-m^-(P_a{\dot C}_aP_a)
+m^+(P_b{\dot C}_bP_b)+\sum_{a<t<b}\sign(P_t{\dot C}_tP_t)\\
&=&-m^+(\gamma(D),\Gamma(\gamma(\ker(A^*+C_s+\delta I)),a)\\
& &+m^-(\gamma(D),\Gamma(\gamma(\ker(A^*+C_s+\delta I)),b)\\
& &-\sum_{a<t<b}\sign(\gamma(D),\Gamma(\gamma(\ker(A^*+C_s+\delta I)),t)\\
&=&-i(\gamma(D),\gamma(\ker(A^*+C_s+\delta I)))\\
&=&-i(\gamma(D),\gamma(\ker(A^*+C_t))),
\ey
where $P_t$ is the orthogonal projection from $H$ onto $\ker(A_D+C_t)$.
\qe

Now we turn to the Maslov-type indices.

\d \lb{d01m3.4}
Let $(X_l,\omega_l)$ be symplectic Hilbert spaces with $\omega_l(x,y)=(j_lx,y)$,
$x,y\in X_l$, $j_l\in\B(X)$ are invertible, and $j_l^*=-j_l$, where $l=1,2$.
Then we have a symplectic Hilbert space
$(H=X_1\oplus X_2,(-\omega_1)\oplus \omega_2)$. Let $W\in{\cal L}(H)$.
Let $M(t)$, $a\le t\le b$ be a curve in $\Sp(X_1,X_2)$ such that
$\Gr(M(t))\in{\cal FL}(W)$ for all $t\in[a,b]$. The {Maslov-type index}
$i_W(M(t))$ is defined to be $i(\Gr(M(t), W)$. If $a=0$, $b=T$, $(X_1,\omega_1)=(X_2,\omega_2)$
and $M(0)=I$, we denote by $\nu_{T,W}(M(t))=\dim_{\C}(\Gr(M(T)\cap W)$.
\ed

The Maslov-type indices have the following property.

\la \lb{l01m3.11}
Let $(X_l,\omega_l)$ be symplectic Hilbert spaces with $\omega_l(x,y)=(j_lx,y)$,
where $x,y\in X_l$, $j_l\in\B(X_l)$ are invertible, and $j_l^*=-j_l$, $l=1,2,3,4$.
Let $W$ be a Lagrange subspace of $(X_1\oplus
X_4,(-\omega_1)\oplus\omega_4)$. Let $\gamma_l\in
C([0,1],\Sp(X_l,X_{l+1}))$, $l=1,2,3$ be syplectic paths such that
$\Gr(\gamma_3(s)\gamma_2(t)\gamma_1(s))\in{\cal FL}(W)$ for all
$(s,t)\in[0,1]\times[0,1]$. Then we have

\ee\lb{eq01m3.21}
i_W(\gamma_3\gamma_2\gamma_1)=i_{W^{'}}(\gamma_2)+i_W(\gamma_3\gamma_2(0)\gamma_1),
\eee
where $W^{'}=\diag(\gamma_1(1),\gamma_3(1)^{-1})W$.
\el

\bp Let $M=\diag(\gamma_1(1),\gamma_3(1)^{-1})$. By the homotopy invariance
rel. endpoints of the Maslov-type indices and Corollary
\ref{c01m3.10}, we have

\y
i_W(\gamma_3\gamma_2\gamma_1)
&=&i_W(\gamma_3(1)\gamma_2\gamma_1(1))+i_W(\gamma_3\gamma_2(0)\gamma_1)\\
&=&i(M\Gr(\gamma_3(1)\gamma_2\gamma_1(1)),M W)+i_W(\gamma_3\gamma_2(0)\gamma_1)\\
&=&i_{W^{'}}(\gamma_2)+i_W(\gamma_3\gamma_2(0)\gamma_1).
\ey\qe

The following properties of fundamental solutions for linear ODE
will be used later.

\la \lb{l01m3.12}
Let $j\in C^1([0,+\infty),\Gl(m,\C))$ be a curve of skew
selfadjoint matrices, and $b\in C([0,+\infty),\gl(m,\C))$ be a
curve of selfadjoint matrices. Let $\gamma\in
C^1([0,+\infty),\Gl(m,\C))$ be the fundamental solution of

\ee \lb{eq01m3.22}
-j{\dot x}-\frac{1}{2}{\dot j}x=bx.
\eee
Then we have $\gamma(t)^*j(t)\gamma(t)=j(0)$ for all $t$.
\el

\bp By the definition of the fundamental solution, we have
$\gamma(0)^*j(0)\gamma(0)=j(0)$. Since $j^*=-j$ and $b^*=b$, we have

\y \frac{d}{dt}(\gamma(t)^*j(t)\gamma(t))
&=&{\dot \gamma }^*j\gamma+\gamma^*{\dot
j}\gamma+\gamma^*j{\dot\gamma}\\
&=&(-b\gamma-\frac{1}{2}{\dot j})^*j^{*-1}j\gamma+\gamma^*{\dot
j}\gamma+\gamma^*jj^{-1}(-b\gamma-\frac{1}{2}{\dot j})\\
&=&\gamma^*(b-\frac{1}{2}{\dot j}+{\dot j}-b-\frac{1}{2}{\dot
j})\gamma\\
&=&0.\ey
So we have $\gamma(t)^*j(t)\gamma(t)=j(0)$.
\qe

\la\lb{l01m3.13} Let $B\in C([0,+\infty),\gl(m,\C))$ and $P\in
C^1([0,+\infty),\Gl(m,\C))$ be two curves of matrices. Let
$\gamma\in C^1([0,+\infty),\Gl(m,\C))$ be the fundamental solution
of

\ee \lb{eq01m3.23}
{\dot x}=Bx,
\eee
and $\gamma^{'}\in C^1([0,+\infty),\Gl(m,\C))$ be the fundamental solution of

\ee \lb{eq01m3.24}
{\dot y}=(PBP^{-1}+{\dot P}P^{-1})y.
\eee
Then we have

\ee \lb{eq01m3.25}
\gamma^{'}=P\gamma P(0)^{-1}.
\eee
\el

\bp Direct calculation shows

$$\frac{d}{dt}(P\gamma P(0)^{-1})=(PBP^{-1}+{\dot P}P^{-1})P\gamma
P(0)^{-1}$$
and $P(0)\gamma P(0)^{-1}=I$. By definition, $P\gamma P(0)^{-1}$
is the fundamental solution of (\ref{eq01m3.24}).
\qe

\c \lb{c01m3.11} Let $j_1,j_2\in C^1([0,+\infty),\Gl(m,\C))$ be two curves of skew
selfadjoint matrices. Let $P\in C^1([0,+\infty),\Gl(m,\C))$ be a curve
of matrices such that $P^*j_2P=j_1$, and $b\in C([0,+\infty),\Gl(m,\C))$
be a curve of selfadjoint matrices.
Let $\gamma\in C^1([0,+\infty),\Gl(m,\C))$ be the fundamental solution
of

\ee \lb{eq01m3.26}
-j_1{\dot x}-\frac{1}{2}{\dot  j}_1x=bx,
\eee
and $\gamma^{'}\in C^1([0,+\infty),\Gl(m,\C))$ be the fundamental solution of

\ee \lb{eq01m3.27} j_2{\dot y}-\frac{1}{2}{\dot
j}_2y=(P^{*-1}bP^{-1}+Q)y, \eee where $Q=\frac{1}{2}(P^{*-1}{\dot
P}^*j_2-j_2{\dot P}P^{-1})$. Then we have

\ee \lb{eq01m3.28}
\gamma^{'}=P\gamma P(0)^{-1}.
\eee

In particular, when $j_1$ and $j_2$ are constant matrices, we have
$$Q=P^{*-1}{\dot P}^*j_2=-j_2{\dot P}P^{-1}.$$
\ec

\bp Take $B=-j_1^{-1}(b+\frac{1}{2}{\dot  j}_1)$ in Lemma \ref{l01m3.13}, we have

\y -j_2(PBP^{-1}+{\dot P}P^{-1})-\frac{1}{2}{\dot  j}_2
&=&-j_2(P(-j_1)^{-1}(b+\frac{1}{2}{\dot  j}_1)P^{-1}+{\dot P}P^{-1})-\frac{1}{2}{\dot
j_2}\\
&=&P^{*-1}(b+\frac{1}{2}{\dot  j}_1)P^{-1}-j_2{\dot
P}P^{-1}-\frac{1}{2}{\dot  j}_2\\
&=&P^{*-1}bP^{-1}-j_2{\dot P}P^{-1}+\frac{1}{2}(P^{*-1}{\dot
j}_1P^{-1}-j_2)\\
&=&P^{*-1}bP^{-1}-j_2{\dot
P}P^{-1}+\frac{1}{2}(P^{*-1}\frac{d}{dt}(P^*j_2P)
P^{-1}-j_2)\\
&=&P^{*-1}bP^{-1}+Q. \ey By Lemma \ref{l01m3.13}, our results
holds.\qe

The following is a special case of the spectral flow formula.

Let $j\in C^1([0,T],\Gl(m,\C))$ be a curve of skew
selfadjoint matrices. Then we have
symplectic Hilbert spaces $(\C^m,\omega(t))$ with standard Hermitian inner product and
$\omega(t)(x,y)=(j(t)x,y)$,
for all $x,y\in \C^m$ and $t\in[0,T]$.
Then we have a symplectic Hilbert space
$(V=\C^m\oplus \C^m,(-\omega(0))\oplus \omega(T))$. Let $W\in{\cal L}(V)$.
Let $b_s(t)\in\B(\C^m)$, $0\le s \le 1$, $0\le t\le T$ be a continuous
family of selfadjoint matrices such that $b_0(t)=0$. By Lemma \ref{l01m3.12}, there are
continuous family of matrices $M_s(t)\in\Gl(m,\C)$ such that $M_s(0)=I$ and
$$-j\frac{d}{d t}M_s(t)-\frac{1}{2}(\frac{d}{d t}j)M_s(t)=b_s(t)M_s(t).$$
Set
\y
X&=&L^2([0,T],\C^m), \quad D_m=H^1_0([0,T],\C^m), \\
D_M&=&H^1([0,T],\C^m),\quad D_W=\{x\in D_M;(x(0),x(t))\in W\}.
\ey
Let $A_M\in{\cal C}(X)$ with domain $D_M$ be defined by
$$A_Mx=-j\frac{d}{d t}x-\frac{1}{2}(\frac{d}{d t}j)x.$$
Set $x\in D_M$, $A=A_M|_{D_m}$, $A_W=A_M|_{D_W}$. Let $C_s\in\B(X)$ be defined
by $(C_sx)(t)=b_s(t)x(t)$, $x\in X$, $t\in[0,T]$.

\p \lb{c01m3.2} Set $W^{'}=\diag(I,M_0(T)^{-1})W$. Then we have

\ee \lb{eq01m3.5}
I(A_W,A_W-C_1)=i_{W^{'}}(M_0^{-1}M_1).
\eee
\ep

\bp The Sobolev embedding theorem shows that $D_M\subset C([0,T],\C^m)$.
For any $x\in D_M$, define $\gamma(x)=(x(0),x(T))$.
Direct calculation shows that $D_M/D_m=\C^m\oplus\C^m$ with
symplectic structure $(\diag(j(0),-j(T))\gamma(x),\gamma(y))$, $x,y\in D_M$,
and $\gamma$ is the abstract trace map. Moreover, $A^*=A_M$,
$\gamma(A^*-C_s)=\Gr(M_s(T))$, and $\gamma(D_W)=W$.
By Proposition \ref{p01m3.1} and Lemma \ref{l01m3.11}, we have

\y I(A_W,A_W-C_1)
&=&-\sf\{A_W-C_s\}\\
&=&i(\{M_s(T);0\le s\le 1\}, W)\\
&=&i_W(M_0(T)(M_0(T)^{-1}M_s(T))I;0\le s\le 1)\\
&=&i_{W^{'}}(M_0(T)^{-1}M_s(T);0\le s\le 1)\\
&=&-i_{W^{'}}(M_0(t)^{-1}M_0(t);0\le t\le T)+i_{W^{'}}(M_0(0)^{-1}M_s(0);0\le s\le 1)\\
& &+i_{W^{'}}(M_0(t)^{-1}M_1(t);0\le t\le T)\\
&=&i_{W^{'}}(M_0^{-1}M_1).
\ey
\qe

\s{Proof of the main results}

In this section we will use the notations in {\S}1.
Let $R$ be a subspace of $\C^{2n}$. Set $H=L^2([0,T],\C^n)$, where $T>0$.
Let $K_R$ be a closed operator on $H$. Its domain is $H_R$ defined by
(\ref{eq01m1.3}), and $H_{R}x=\dot x$ for all $x\in H_R$.
Set
$$ X=L^2([0,T],\C^{2n}), \quad
D_{W(R)}=\{x\in H^1([0,T],\C^{2n});(x(0),x(t))\in W(R)\}.
$$
For any $x,y\in H^1([0,T],\C^{2n})$, define
$$(x,y)_1=(x,y)+(\dot x,\dot y).$$
Let $A_{W(R)}\in{\cal C}(X)$ with domain $D_{W(R)}$ be defined by
$A_{W(R)}x=-J\frac{d}{d t}x$, $x\in D_{W(R)}$. Let $C_s\in\B(X)$ be defined
by $(C_sx)(t)=b_s(t)x(t)$, $x\in X$, $t\in[0,T]$. Then we have $K_R^*=-K_{R^b}$.
Consider the standard orthogonal decomposition
$$\C^{2n}=(\C^n\times\{0\})\oplus(\{0\}\times \C^n).$$
It induces orthogonal decompositions $X=H\oplus H$ and
$D_{W(R)}=H_{R^b}\oplus H_{R}$. Under such orthogonal decompositions,
$A_{W(R)}$ is in block form $A_{W(R)}=\pmatrix{0&K_R\cr K_R^*&0\cr}$.
Let $b_s(t), C_s$ be in block form
\y
b_s(t)&=&\pmatrix{b_{11}(s,t)&b_{12}(s,t)\cr b_{21}(s,t)&b_{22}(s,t)\cr}\\
C_s&=&\pmatrix{C_{11}(s)&C_{12}(s)\cr C_{21}(s)&C_{22}(s)\cr}.
\ey
Define $P_s,Q_s,R_s\in\B(H_R)$ by
\y
(P_sx,y)_1&=&(C_{11}(s)^{-1}K_Rx,K_Ry),\\
(Q_sx,y)_1&=&-(C_{11}(s)^{-1}C_{12}(s)x,K_Ry),\\
(R_sx,y)_1&=&((C_{21}(s)C_{11}(s)^{-1}C_{12}(s)-C_{22}(s))x,y)
\ey
for all $x,y\in\H_R$.
Then the index forms ${\cal I}_s$ defined by (\ref{eq01m1.4})
satisfy

$${\cal I}_s(x,y)=((P_s+Q_s+Q_s^{*}+R_s)x,y)_1.$$

\la \lb{l01m4.1}
The operator $P_s+Q_s+Q_s^{*}+R_s\in\B(H_R)$ is a curve of
Fredholm selfadjoint operators.
\el

\bp Since ${\cal I}_s$ are bounded symmetric quadratic forms on $H_R$, by Riesz
representation theorem, $P_s+Q_s+Q_s^{*}+R_s\in\B(H_R)$. Similarly we can see
that they form a continuous curve.

Let $Q_s^{'}=-C_{11}(s)^{-1}C_{12}(s)$ and
$R_s^{'}=C_{21}(s)C_{11}(s)^{-1}C_{12}(s)-C_{22}(s)$ be two bounded operators
on $H$. Then $Q_s^*=(K_R^*K_R+I)^{-1}Q_s^{'*}K_R$ and
$R_s=(K_R^*K_R+I)^{-1}R_s^{'}$. So $Q_s^*$ and $R_s$ maps bounded subset of
$H_R$ into the bounded subset in the domain $\Da(K_R^*K_R)$ of $K_R^*K_R$.
Since $\Da(K_R^*K_R)$ is a closed subspace of $H^2([0,T],\C^n)$, by
Sobolev embedding theorem, the embedding of $\Da(K_R^*K_R)$ into $H_R$ is compact.
So $Q_s^*$, $Q_s$, and $R_s$ are compact.

Now we prove that $P_s$ is Fredholm and then our lemma is proved.
If $p_s(t)$ is positive definite,
we can choose $q_s(t)=0$ and a positive definite curve $r_s(t)$
such that $P_s+R_s$ is positive definite. So $P_s+R_s$ is invertible.
Since $R_s$ is compact, $P_s$ is Fredholm.
Here it is only required that
$p_s(t)$ continuous. In the general case, we have to assume that
$p_s(t)$ is $C^1$ in $t$. Note that $H_{\C^{2n}}=H^1([0,T],\C^n)$.
Consider the operator $p_s:H_{\C^{2n}}\to H_{\C^{2n}}$. Let
$j:H_R\to H_{\C^{2n}}$ be the injection.
Then $p_s$ is invertible and $p_sj$ is Fredholm. For any $x\in H_R$ and
$y=H_{\C^{2n}}$, the inner product $((P_s-p_s)x,y)_1$ consists only the
lower-order terms (i.e., no second-order differential involved)
and some boundary terms. Similar to the above proof, we can conclude that the
lower-order terms correspond to compact operators. The boundary terms correspond
to finite rank operators. So $jP_s-p_sj$ is compact. Since $p_sj$ and $j$ are
Fredholm, $jP_s$ and $P_s$ are Fredholm.
\qe

The following lemma is the key to the proof of Theorem \ref{th01m1.1}.

\la \lb{l01m4.2} Let $u_{b_s}(x)=(p_sK_Rx+q_sx,x)$ for all $x\in H_R$ and
$0\le s\le 1$. Then we have

\ee \lb{eq01m4.1}
\ker(A_{W(R)}-C_s)=\{u_{b_s}(x);x\in\ker\;{\cal I}_s\}.
\eee
Moreover, for any $x,y\in H_R$, we have

\ee \lb{eq01m4.2}
-\left(\left(\frac{d}{ds}C_s\right)u_{b_s}(x), u_{b_s}(y)\right)
= \left(\left(\frac{d}{ds}{\cal I}_s\right)x,y\right).
\eee
\el
\qe

\bp Since ${\cal I}_s(x,y)=(p_sK_Rx+q_sx,K_Ry)+(q_s^*K_Rx+r_sx,y)$
for all $x,y\in H_R$, by the definition of $K_R^*$, we have $x\in\ker\;{\cal
I}_s$ if and only if $x\in\Da(K_R^*)$, and
$K_R^*(p_sK_Rx+q_sx)+(q_s^*K_Rx+r_sx)=0$, if and only if
$u_{b_s}(x)\in\ker(A_{W(R)}-C_s)$. So equation (\ref{eq01m4.1}) is
proved.

Now we turn to equation (\ref{eq01m4.2}). Set
$Z_s=\pmatrix{p_s&q_s\cr0&1\cr}$. By the definition of $b_s$ and direct computation we
have

\ee\lb{eq01m4.2a}
-Z_s^*b_sZ_s=\pmatrix{-p_s&0\cr 0&r_s\cr},\qquad
-Z_s^*b_s=\pmatrix{-I&q_s\cr 0&r_s\cr},\qquad
-b_sZ_s=\pmatrix{-I&0\cr  q_s^*&r_s\cr}.\eee
So we have

\y
-Z_s^*\left(\frac{d}{ds}b_s\right)Z_s
&=&-\frac{d}{ds}\left(Z_s^*b_sZ_s\right)
+Z_s^*b_s\frac{d}{ds}Z_s+Z_s^*\frac{d}{ds}(b_sZ_s)\\
&=&\pmatrix{-\frac{d}{ds}p_s&0\cr 0&\frac{d}{ds}r_s\cr}
-\pmatrix{-I&q_s\cr 0&r_s\cr}\pmatrix{\frac{d}{ds}p_s&\frac{d}{ds}q_s\cr 0&0\cr}
-\pmatrix{\frac{d}{ds}p_s&0\cr \frac{d}{ds}q_s^*&0\cr}\pmatrix{-I&0\cr
q_s^*&r_s\cr}\\
&=&\pmatrix{-\frac{d}{ds}p_s&0\cr 0&\frac{d}{ds}r_s\cr}
-\pmatrix{-\frac{d}{ds}p_s&-\frac{d}{ds}q_s\cr 0&0\cr}
-\pmatrix{-\frac{d}{ds}p_s&0\cr -\frac{d}{ds}q_s^*&0\cr}\\
&=&\frac{d}{ds}\pmatrix{p_s&q_s\cr q_s^*&r_s\cr}.\ey
Hence for all $x,y\in H_R$, we have

\y
-\left(\left(\frac{d}{ds}C_s\right)u_{b_s}(x), u_{b_s}(y)\right)
&=&-\int_0^T \left(\left(\frac{d}{ds}b_s\right)Z_s\pmatrix{{\dot x}\cr
x\cr},Z_s\pmatrix{{\dot x}\cr x\cr}\right)dt\\
&=&\int_0^T \left(\frac{d}{ds}\pmatrix{p_s&q_s\cr q_s^*&r_s\cr}\pmatrix{{\dot x}\cr
x\cr},\pmatrix{{\dot x}\cr x\cr}\right)dt\\
&=& \left(\left(\frac{d}{ds}{\cal I}_s\right)x,y\right).
\ey
\qe

Now we can prove Theorem \ref{th01m1.1}.

{\bf Proof of Theorem \ref{th01m1.1}.}\hspace{2mm}
By Lemma \ref{l01m4.1}, $\sf\{{\cal I}_s\}$ is well-defined. Since
the spectral flow is invariant under homotopy with endpoints fixed, by small
perturbation which fixes the endpoint, we can assume that
$C_s$ is a $C^1$ curve. By Theorem 4.22 in \cite{RS2}, we can take a sufficient
small $\delta>0$ such that $C_{11}(s)-\delta I$ is invertible,
$A_{W(R)}-C_s+\delta I$, $0\le s\le 1$ has only regular crossing,
and $A_{W(R)}-C_s+aI$ is invertible for $s=0,1$ and $a\in(0,\delta]$.
Let ${\cal I}_{s,a}$ be the correspondent index form of $C_s-a I$.
By Proposition \ref{p01m2.1}, Proposition \ref{c01m3.2} and Lemma \ref{l01m4.2}
we have

\y I({\cal I}_0,{\cal I}_1)
&=&-\sf\{{\cal I}_s,0\le s\le 1\} \\
&=&-\sf\{{\cal I}_{s,\delta},0\le s\le 1\}\\
&=&-\sf\{A_{W(R)}-C_s+\delta I,0\le s\le 1\}\\
&=&-\sf\{A_{W(R)}-C_s,0\le s\le 1\}\\
&=&i_{W(R)}(\gamma_1)-i_{W(R)}(\gamma_0).
\ey
\qe

To prove Theorem \ref{th01m1.2}, we need some preparations.

\la \lb{l01m4.4}
Theorem \ref{th01m1.2} holds in the following two cases:
\begin{description}
\item[(i)] $R=\C^{2n}$,
\item[(ii)] $p\in C([0,T],\gl(n,\C))$ is a path of positive definite matrices,
and $P(t)=\int_0^tp(s)d s$ for all $t\in[0,T]$.
\end{description}
\el

\bp By the definition of $W(R)$ we have
$$\Gr(\gamma(t))\cap W(R)=\{(x,y,x,y+P(t)x);(x,x)\in R^b, (y,y+P(t)x)\in R\}.$$
In the case (i), we have $R^b=\{0\}$ and
$\Gr(\gamma(t))\cap W(R)=\{(0,y,0,y);y\in\C^n\}$ for all
$t\in[0,T]$. So we have $i_{W(R)}(\gamma)=0$ and Theorem
\ref{th01m1.2} holds.

Now we consider the case (ii). Since $P(t)$ is positive definite for all
$t\in(0,T]$, for all $(x,y,x,y+P(t)x)\in\Gr(\gamma(t))\cap W(R)$,
we have $-(x,P(t)x)=((x,-x),(y,y+P(t)x))=0$ and hence $x=0$. So we
have
$$\Gr(\gamma(t))\cap W(R)=\cases{\Gr(I)\cap W(R), &if $t=0$,\cr
\{(0,y,0,y);(y,y)\in R\}, &if $t\in(0,T]$.}$$
By Corollary \ref{c01m3.13} we have
\y i_{W(R)}(\gamma)&=&\dim_{\C}(\Gr(\gamma(0))\cap W(R))
-\dim_{\C}(\Gr(\gamma(T))\cap W(R))\\
&=&\dim_{\C}(\Gr(I)\cap W(R))-\dim_{\C}(\Gr(I)\cap R)\\
&=&\dim_{\C}(\Gr(I)\cap R^b)\\
&=&\dim_{\C}S.\ey
\qe

\p \lb{p01m4.1} Let $\gamma\in C([0,T],\Sp(2n,\C))$ be such that
$\gamma(0)=I$. Let $R_1\subset R_2$ be two linear subspaces of $\C^{2n}$. Define
$${\cal N}=\{(x,y,z,u)\in \Gr(\gamma(T));
(x,z)\in R_1^b,(y,u)\in R_2\},$$
and
$$Q((x_1,y_1,z_1,u_1),(x_2,y_2,z_2,u_2))=(z_1,u_2)-(x_1,y_2)$$
for all $(x_1,y_1,z_1,u_1),(x_2,y_2,z_2,u_2)\in{\cal N}$. Then
$Q$ is a quadratic form on ${\cal N}$, and we have
\ee \lb{eq01m4.30}
i_{W(R_2)}(\gamma)-i_{W(R_1)}(\gamma)=C(\gamma(T);R_1,R_2)+
\dim_{\C}(\Gr(I)\cap R_2^b)-\dim_{\C}(\Gr(I)\cap R_1^b),
\eee
where
$$C(\gamma(T);R_1,R_2)=m^-(Q)+\dim_{\C}\ker Q-\dim_{\C}(\Gr(\gamma(T))\cap
W(R_2)).$$
We call $C(\gamma(T);R_1,R_2)$ the {\bf Morse concavity} of $\gamma(T)$ with respect
to $R_1,R_2$.\ep

\bp By \cite{Wi}, there exist paths of matrices $p,q,r\in\gl(n,\C)$ such
that $p(t)$ are positive definite, $r(t)=r(t)^*$ for all
$t\in[0,T]$, and $\gamma_1(T)=\gamma(T)$ if $p_s=p$, $q_s=sq$,
$r_s=sr$ for all $s\in[0,1]$, and $\gamma_s$ are the fundamental
solution of (\ref{eq01m1.2}).
Let ${\cal I}_s$ be defined by (\ref{eq01m1.4}). By Theorem
\ref{th01m1.1} and Lemma \ref{l01m4.4}, we have

\ee \lb{eq01m4.31}
m^-({\cal I}_1|_{H_{R_k}})=i_{W(R_k)}(\gamma_1)-\dim_{\C}(\Gr(I)\cap R_k^b),
\quad k=1,2.
\eee
Let $N$ be the ${\cal I}_1$-complement of $H_{R_1}$ in $H_{R_2}$.
Direct computation shows that

$$N=\left\{
\begin{array}{c}x\in H_{R_2};(u_{b_1}(x))(t)=\gamma_1(t)(u_{b_1}(x))(0)\;
\ox{\rm for}\;\ox{\rm all}\;t\in[0,T]\\
\ox{\rm and}\;(p(0){\dot x}(0)+q(0)x(0),p(T){\dot x}(T)+q(T)x(T))\in R_1^b
\end{array}\right\}.$$
Define $\varphi:N\to {\cal N}$ by
$$\varphi(x)=((u_{b_1}(x))(0),(u_{b_1}(x))(T)).$$
Then $\varphi$ is a linear isomorphism. Direct computation shows
that
$$ {\cal I}_1(x,y)=Q(\varphi(x),\varphi(y))$$
for all $x,y\in N$. By Proposition \ref{p01m2.2} and Lemma \ref{l01m4.2}
we have
\aa
m^-({\cal I}_1|_{H_{R_2}})-m^-({\cal I}_1|_{H_{R_1}})
&=&m^-({\cal I}_1|_N)+\dim_{\C}\ker({\cal I}_1|_N)-\dim_{\C}\ker({\cal
I}_1|_{H_{R_2}})\nn\\
&=&m^-(Q)+\dim_{\C}\ker Q-\dim_{\C}(\Gr(\gamma_1(T))\cap W(R_2))\nn\\
&=&m^-(Q)+\dim_{\C}\ker Q-\dim_{\C}(\Gr(\gamma(T))\cap W(R_2))\nn\\
&=&C(\gamma(T);R_1,R_2).
\lb{eq01m4.32} \eaa
By the fact that $\gamma$ and $\gamma_1$ has the same end points,
we have
$$i_{W(R_1)}(\gamma)-i_{W(R_1)}(\gamma_1)=i_{W(R_2)}(\gamma)-i_{W(R_2)}(\gamma_1).$$
By (\ref{eq01m4.31}) and (\ref{eq01m4.32}), we have
\y i_{W(R_2)}(\gamma)-i_{W(R_1)}(\gamma)
&=&i_{W(R_2)}(\gamma_1)-i_{W(R_1)}(\gamma_1)\\
&=&(m^-({\cal I}_1|_{H_{R_2}})+\dim_{\C}(\Gr(I)\cap R_2^b))\\
& &-(m^-({\cal I}_1|_{H_{R_1}})+\dim_{\C}(\Gr(I)\cap R_1^b))\\
&=&C(\gamma(T);R_1,R_2)+\dim_{\C}(\Gr(I)\cap R_2^b)-\dim_{\C}(\Gr(I)\cap R_1^b).
\ey
\qe

{\bf Proof of Theorem \ref{th01m1.2}.}\hspace{2mm} Firstly we assume that
$P(0)=0$. Set $R_1=R$ and $R_2=\C^{2n}$. Let $Q$, ${\cal N}$ be defined by Proposition
\ref{p01m4.1}. By the definition of ${\cal N}$ we have
$${\cal N}=\{(x,y,x,P(T)x+y);(x,x)\in R^b\}.$$
Define $\varphi:S\times \C^n\to{\cal N}$ by
$$\varphi(x,y)=(x,y,x,P(T)x+y).$$
Then $\varphi$ is a linear isomorphism. So we have $\dim_{\C}{\cal
N}=\dim_{\C}S+n$.
By the definition of $Q$ we have
\y Q(\varphi(x_1,y_1),\varphi(x_2,y_2))
&=&(x_1,P(T)x_2+y_2)-(x_1,y_2)\\
&=&(P(T)x_1,x_2).
\ey
So we have $m^+(Q)=m^+(P(T)|_S)$. By the definition of
$C(\gamma(T);R,\C^{2n})$ we have
\y
C(\gamma(T);R,\C^{2n})&=&m^-(Q)+\dim_{\C}\ker Q-\dim_{\C}(\Gr(\gamma(T))\cap
W(R))\\
&=&\dim_{\C}{\cal N}-m^+(Q)-n\\
&=&\dim_{\C}S-m^+(P(T)|_S).
\ey
By Proposition \ref{p01m4.1} we have
\y i_{W(R)}(\gamma)
&=&i_{W(\C^{2n})}(\gamma)-C(\gamma(T);R,\C^{2n})
-\dim_{\C}(\Gr(I)\cap (\C^{2n})^b)+\dim_{\C}(\Gr(I)\cap R^b)\\
&=&0-(\dim_{\C}S-m^+(P(T)|_S))-0+\dim_{\C}S\\
&=&m^+(P(T)|_S).\ey

Now we consider the general case. Define
$M_t(s)=\pmatrix{I&0\cr s P(t)&0\cr}$
for all $s\in[0,1]$ and $t\in[0,T]$.
Then $M_t(s)\in\Sp(2n,\C)$, and we have
\y i_{W(R)}(\gamma)
&=&i_{W(R)}(M_T)-i_{W(R)}(M_0)\\
&=&m^+(P(T)|_S)-m^+(P(0)|_S).
\ey
\qe

Now we prove Theorem \ref{th01m1.3}. Let $a$, $p_1^{'}$, $q_1^{'}$, $r_1^{'}$,
$b_1^{'}$, $R^{'}$ be as in {\S}1. Let $C_1^{'}$ be the corresponding
bounded operator of $C_1$. The following lemma follows from direct calculation.

\la \lb{l01m4.3} We have

\aa
\lb{eq01m4.3}
\pmatrix{p_1^{'}&q_1^{'}\cr q_1^{'*}&r_1^{'}\cr}
&=&\pmatrix{a^{*}&0\cr {\dot a}^{*}&a^{*}\cr}
\pmatrix{p_1&q_1\cr q_1^{*}&r_1\cr}
\pmatrix{a&\dot a\cr 0&a\cr},\\
\lb{eq01m4.4} b_1^{'}&=&\diag(a^{-1},a^{*})b_1\diag(a^{*-1},a)
+\pmatrix{0&-a^{-1}{\dot a}\cr -{\dot a}^{*}a^{*-1}&0\cr},\\
\lb{eq01m4.5}
A_{W(R^{'})}-C_1^{'}&=&\diag(a^{-1},a^{*})(A_{W(R)}-C_1)\diag(a^{*-1},a),\\
\lb{eq01m4.7}
a^{-1}K_Ra&=&K_{R^{'}}+a^{-1}{\dot a}.
\eaa
\el\qe

By Corollary \ref{c01m3.11} we have

\c\lb{c01m4.10} We have \ee \lb{eq01m4.20}
\gamma_1^{'}=\diag(a^*,a^{-1})\gamma_1\diag(a(0)^{*-1},a(0)). \eee
\ec






\la \lb{l01m4.5} Let $a\in C([0,T],\Gl(n,\C))$. Set
$\gamma=\diag(a^*,a^{-1})$. Then $\gamma\in C([0,T],\Sp(2n,\C))$, and we have
\aa
i_{W(R)}(\gamma)
&=&\dim_{\C}(\Gr(a(0)^{-1})\cap R)-\dim_{\C}(\Gr(a(T)^{-1})\cap R)\nn\\
&=&\dim_{\C}(\Gr(a(0)^*)\cap R^b)-\dim_{\C}(\Gr(a(T)^{*})\cap R^b).
\lb{eq01m4.34}\eaa
\el

\bp Clearly we have $\gamma\in C([0,T],\Sp(2n,\C))$. We divide the
proof into three steps.

{\bf Step 1.} $a(0)=I$ and $a\in C^1([0,T],\Gl(n,\C))$.

In this case, set
$$c=\pmatrix{0&-a^{-1}{\dot a}\cr -{\dot a}^*a^{*-1}&0\cr}.$$
Define $C\in\B(X)$ by $(C u)(t)=c(t)u(t)$. Then $\gamma$, $a^{-1}$
and $a^*$ are the fundamental solutions of ${\dot u}=J c u$,
${\dot x}=-a^{-1}{\dot a}x$ and ${\dot x}={\dot a}^*a^{*-1}x$
respectively. Since $K_R^*=-K_{R^b}$, we have \y
\dim_{\C}\ker(K_R+a^{-1}{\dot a})&=&\dim_{\C}(\Gr(a(T)^{-1})\cap
R),\\
\dim_{\C}\ker(K_R^*+{\dot
a}^*a^{*-1})&=&\dim_{\C}(\Gr(a(T)^{*})\cap R^b). \ey By
Proposition \ref{c01m3.2} and Lemma \ref{l01m2.8}, we have
\y i_{W(R)}(\gamma)&=&I(A_{W(R)},A_{W(R)}-C)\\
&=&-\sf\{A_{W(R)}-s C;0\le s\le 1\}\\
&=&-\sf\{\pmatrix{0&-K_R-sa^{-1}{\dot a}\cr -K_R^*-s{\dot
a}^*a^{*-1}&0\cr};0\le s\le 1\}\\
&=&\dim_{\C}\ker(K_R)-\dim_{\C}\ker(K_R+a^{-1}{\dot a})\\
&=&\dim_{\C}\ker(K_R^*)-\dim_{\C}\ker(K_R^*+{\dot a}^*a^{*-1})\\
&=&\dim_{\C}(\Gr(I)\cap R)-\dim_{\C}(\Gr(a(T)^{-1})\cap R)\\
&=&\dim_{\C}(\Gr(I)\cap R^b)-\dim_{\C}(\Gr(a(T)^{*})\cap R^b).
\ey

{\bf Step 2.} $a(0)=I$.

Since $\Gl(n,\C)$ is a connected Lie group, there exists
$H_s(t)\in\Gl(n,\C)$ such that $H_0(t)=a(t)$, $H_1$ is smooth, $H_s(0)=I$, and
$H_s(T)=a(T)$ for all $s\in[0,1]$, $t\in[0,T]$. By Step 1 we have
\y i_{W(R)}(\gamma)&=&i_{W(R)}(\diag(H_1^*,H_1^{-1}))\\
&=&\dim_{\C}(\Gr(I)\cap R)-\dim_{\C}(\Gr(a(T)^{-1})\cap R)\\
&=&\dim_{\C}(\Gr(I)\cap R^b)-\dim_{\C}(\Gr(a(T)^{*})\cap R^b).
\ey

{\bf Step 3.} General case.

Since $\Gl(n,\C)$ is a connected Lie group, there exists
$\alpha\in C([0,T],\Gl(n,\C))$ such that $\alpha(0)=I$,
$\alpha(T)=a(0)$.
By Step 2 we have
\y i_{W(R)}(\gamma)
&=&(i_{W(R)}(\diag(\alpha^*,\alpha^{-1}))+i_{W(R)}(\gamma))
-i_{W(R)}(\diag(\alpha^*,\alpha^{-1}))\\
&=&(\dim_{\C}(\Gr(I)\cap R)-\dim_{\C}(\Gr(a(T)^{-1})\cap R))\\
& &-(\dim_{\C}(\Gr(I)\cap R)-\dim_{\C}(\Gr(a(0)^{-1})\cap R))\\
&=&\dim_{\C}(\Gr(a(0)^{-1})\cap R)-\dim_{\C}(\Gr(a(T)^{-1})\cap R)\\
&=&\dim_{\C}(\Gr(a(0))\cap R^b)-\dim_{\C}(\Gr(a(T)^{*})\cap R^b).
\ey
\qe

{\bf Proof of Theorem \ref{th01m1.3}.}\hspace{2mm}
By Corollary \ref{c01m4.10}, Lemma \ref{l01m3.11} and Lemma \ref{l01m4.5},
we have
\y i_{W(R^{'})}(\gamma_1^{'})
&=&
i_{W(R^{'})}(\diag(a^*,a^{-1})\gamma_1\diag(a(0)^{*-1},a(0)))\\
&=&i_{W(R)}(\gamma_1)+i_{W(R^{'})}(\diag(a^*,a^{-1})\diag(a(0)^{*-1},a(0))^{-1})\\
&=&i_{W(R)}(\gamma_1)+\dim_{\C}(\Gr(a(0)^{-1}a(0))\cap R^{'})
-\dim_{\C}(\Gr(a(T)^{-1}a(0))\cap R^{'})\\
&=&i_{W(R)}(\gamma_1)+\dim_{\C}(\Gr(I)\cap R^{'})-\dim_{\C}(\Gr(I)\cap
R).
\ey
\qe





\bibliographystyle{abbrv}

\begin{thebibliography}{99}


\bibitem{AgSa} A. A. Agrachev, A. V. Sarychev, Abnormal sub-Riemannian geodesics:
Morse index and rigidity. {\it Ann. Inst. Henri Poincar\'e,
Analyse non lineair\'e}.
13 (1996). 635-690.

\bibitem{Am} W. Ambrose, The index theorem in Riemannian geometry.
{\it Ann. of Math.} 73 (1961). 49-86.

\bibitem{A} V.I. Arnol'd, Characteristic class entering quantization
conditions. {\it Funkts. Anal. Priloch.} 1 (1967). 1-14 (Russian).
Funct. Anal. Appl. 1 (1967). 1-13 (English transl.).

\bibitem{APS1} M. F. Atiyah, V. K. Patodi, and I. M. Singer,  Spectral
asymmetry and Riemannian geometry. III. {\it Proc. Camb. Phic. Soc.} 79
(1976). 71-99.

\bibitem{BoFu98} B. Booss and K. Furutani, The
Maslov index -- a functional analytical definition and the
spectral flow formula, {\it Tokyo J. Math.}  21 (1998).
1--34.

\bibitem{CLM} S. E. Cappell, R. Lee, and E. Y. Miller,  On the Maslov
index. {\it Comm. Pure Appl. Math.} 47. (1994). 121-186.

\bibitem{DZ1} X. Dai and W. Zhang, Splitting of the familly index. {\it Comm.
Math. Phys.} 182 (1996). 303-317.

\bibitem{DZ2} X. Dai and W. Zhang,  Higher spectral flow.
{\it J. Funct. Analysis.} 157 (1998). 432-469.

\bibitem{Du} J. J. Duistermaat, {\it On the Morse index in variational
calculus.} Adv. Math. 21. (1976). 173-195.

\bibitem{Ho} L. H{\"o}rmander, Fourier integral Operators I. {\it Acta Math.}
127(1971). 79-183.

\bibitem{Ka} T. Kato, Pertubation Theory for Linear Operators.
Springer-Verlag. Berlin. 1980.

\bibitem{Lo7} Y. Long, Bott formula of the Maslov-type index theory.
{\it Pacific J. Math.} 187 (1999). 113-149.

\bibitem{MeP} R. B. Melrose and P. Piazza, Families of Dirac operators,
boundaries and the $b$-calculus. {\it J. Diff. Geom.} 46 (1997). 99-180.

\bibitem{LZh} Y. Long and C. Zhu, Maslov-type index theory for
symplectic paths and spectral flow (II). {\it Chinese Ann. of
Math.} 21B:1 (2000). 89-108.

\bibitem{Mo} M. Morse, The Calculus of Variations in the Large.
A.M.S. Coll. Publ., Vol.18, Amer. Math. Soc., New York, 1934.

\bibitem{PiT1} P. Piccione and D. V. Tausk, The Maslov index and a
generalized Morse index theorem for non-positive definite metrics.
{\it C. R. Acad. Sci. Paris S{\'e}r. I Math.} 331(2000). 385-389.

\bibitem{PiT2} P. Piccione and D. V. Tausk, The Morse index theorem
in semi-Riemannian Geometry. http://xxx.lanl.gov, math.DG/0011090.
{\it Topology.} To appear.

\bibitem{PiT2} P. Piccione and D. V. Tausk, An index theory for
paths that are solutions of a class of strongly indefinite
variational problems. http://xxx.lanl.gov, math.DG/0108044 v1.

\bibitem{RS2} J. Robbin and D. Salamon, The spectral flow
and the Maslov index.
{\it Bull. London Math. Soc.} 27 (1995). 1--33.

\bibitem{Sm} S. Smale, On the Morse index theorem. {\it J. Math. Mech.}
14(1965). 1049-1056.

\bibitem{Uh} K. Uhlenbeck, The Morse index theorem in Hilbert space.
{\it J. Diff.Geom.} 8 (1973). 555-564.

\bibitem{Wi} B. Wilking, Index parity of closed geodesics and
rigidty of Hopf fibrations. {\it Invent. math.} 2001.
DOI 10.1007/s002220100123.

\bibitem{Zh} C. Zhu, Maslov-type index theory and closed characteristic
on compact convex hypersurfaces in $\R^{2n}$. Ph. D. Thesis. Nankai Institue
of Mathematics.

\bibitem{ZhL} C. Zhu and Y. Long, Maslov-type index theory for
symplectic paths and spectral flow (I). {\it Chinese Ann. of Math.} 20B:4
(1999).
413-424.


\end{thebibliography}

\end{document}